\documentclass[11pt]{article}

\usepackage[numbers,sort&compress]{natbib}
\usepackage{graphicx}
\usepackage{amsmath,amssymb,amsfonts}
\usepackage{amsthm}
\usepackage{mathtools}
\usepackage{tikz}
\usepackage{pgfplots}
\pgfplotsset{compat=1.17}
\usepackage{float}
\usepackage{multirow}
\usepackage{booktabs}
\usepackage{url}
\usepackage{enumitem}
\usepackage{blkarray}
\usepackage{rotating}
\usepackage{pdflscape}
\usepackage{hyperref}
\usepackage[margin=1in]{geometry}
\usepackage{authblk}

\theoremstyle{plain}
\newtheorem{lem}{Lemma}
\newtheorem{thm}{Theorem}
\newtheorem{cor}{Corollary}

\theoremstyle{definition}

\newtheorem{assumption}{Assumption}
\newtheorem{ex}{Example}
\newtheorem{remark}{Remark}

\newcommand{\conv}{\mathop{conv}}

\newcommand{\intr}{\mathop{\text{int}}}
\newcommand{\relint}{\mathop{\text{relint}}}

\newcommand{\D}{\mathcal{D}}
\newcommand{\Ds}{\mathcal{D}^*}
\newcommand{\DD}{\mathcal{DD}}
\newcommand{\DDs}{\mathcal{DD}^*}
\newcommand{\K}{\mathcal{K}}

\newcommand{\rank}{\text{rank }}
\newcommand{\Srp}{\mathbb{S}_{+}^{r}}
\newcommand{\Snop}{\mathbb{S}_{+}^{n+1}}
\newcommand{\Srop}{\mathbb{S}_{+}^{r+1}}
\newcommand{\Sno}{\mathbb{S}^{n+1}}
\newcommand{\norm}[1]{\left\|#1\right\|}

\newcommand{\Sn}{{\mathbb{S}^{n}}}
\newcommand{\Snp}{{\mathbb{S}_{+}^{n}}}
\newcommand{\Rn}{\mathbb{R}^{n}}
\newcommand{\Rm}{\mathbb{R}^{m}}

\newcommand{\R}{\mathbb{R}}

\newcommand{\vs}{\mid}
\DeclareMathOperator{\arrow}{arrow}

\newcommand{\aff}{\text{aff }}

\newcommand{\Diag}{{\text{Diag}}}

\newcommand{\ones}{\textbf{1}}
\newcommand{\zeros}{\textbf{0}}
\newcommand{\SEP}{SEP}
\newcommand{\linelabel}[1]{}
\let\cline\cmidrule

\title{Affine Facial Reduction for Semidefinite Relaxations of Binary and Mixed-Binary Optimization Problems}
\author[1]{Hao Hu}
\author[1]{Boshi Yang}
\affil[1]{School of Mathematical and Statistical Sciences, Clemson University, Clemson, SC, USA\\
\texttt{hhu2@clemson.edu}; \texttt{boshiy@clemson.edu}}
\date{}

\hypersetup{
  pdftitle={Affine Facial Reduction for Semidefinite Relaxations of Binary and Mixed-Binary Optimization Problems},
  pdfauthor={Hao Hu and Boshi Yang}
}

\begin{document}

\maketitle

\begin{abstract}
Semidefinite programming (SDP) relaxations can provide strong bounds for binary and mixed-binary optimization problems, but their practical use is often limited by matrix variables of large order and by failures of Slater's condition, which can contribute to numerical difficulties for interior-point solvers. We propose \emph{affine facial reduction (affine FR)}, an automatic, LP-based preprocessing method for SDP relaxations of such problems. The method uses the affine hull of the linear programming relaxation to construct a positive semidefinite face containing the lifted feasible set and an associated facial range vector, thereby replacing the original positive semidefinite matrix variable by one of lower order. The required affine-hull information is obtained by solving a single linear programming problem, followed by linear-algebraic computations. We analyze the relation between affine FR, analytical facial reduction, the Permenter--Parrilo partial FR method, and Sieve-SDP, and establish explicit comparisons among the resulting matrix orders for the Shor relaxation. On MIPLIB 2017 mixed-binary instances, affine FR reduces matrix order more often and by a larger amount than the partial-FR variants considered, with comparable preprocessing times. On selected SAT relaxations with substantial reductions, affine FR can shorten SDP solution times and yield more reliable solver outcomes when the unreduced formulations encounter numerical difficulties.
\end{abstract}



\section{Introduction}
\label{intro}
Binary and mixed-binary optimization models encode yes/no decisions and arise in applications such as radio-frequency assignment, timetabling, scheduling, transportation, network design, and resource allocation. Solving large-scale instances of these problems can be difficult, motivating the development of strong convex relaxations.

Semidefinite programming (SDP) provides a general framework for constructing such relaxations. An SDP has a linear objective function and affine constraints over the cone of positive semidefinite matrices; see \cite{anjos2011handbook,wolkowicz2012handbook}. SDP has attracted significant attention in combinatorial optimization over the past three decades because SDP relaxations can provide strong bounds on optimal values. Notable success has been achieved using SDP relaxations for challenging problems such as the max-cut problem \cite{goemans1995improved,rendl2010solving} and the quadratic assignment problem \cite{zhao1998semidefinite}. \linelabel{newref}\label{letscheck}A specialized first-order method for optimizing lift-and-project relaxations of binary integer programs is proposed in \cite{burer2006solving}. Despite recent advances in first-order methods that have improved the practicality of SDP-based approaches, significant challenges remain, particularly in handling large-scale problems and addressing numerical difficulties associated with failures of regularity conditions \cite{gaar2024strong,krislock2017biqcrunch,oliveira2018admm,zheng2017fast}.
More precisely, the order of the lifted matrix variable grows with the number of variables, and it is well known that SDP solvers do not scale as effectively as linear programming (LP) solvers. Due to the inherent structure of many binary and mixed-binary formulations, the corresponding SDP relaxations often fail to satisfy regularity conditions, such as Slater's condition, which can lead to numerical difficulties.

Facial reduction is a preprocessing technique introduced by Borwein and Wolkowicz to regularize SDPs by restoring Slater's condition \cite{borwein1981regularizing,borwein1981facial}.  It has been applied in several structured settings, including the quadratic assignment problem \cite{zhao1998semidefinite}. In exact arithmetic, facial reduction identifies the minimal face containing the feasible region and restores strict feasibility relative to that face, but computing suitable exposing vectors can be difficult, prompting the development of specialized methods \cite{permenter2018partial,zhu2019sieve}.

In this paper, we introduce \emph{affine FR}, an automatic, LP-based preprocessing method tailored for SDP relaxations of binary and mixed-binary optimization problems. Our approach uses affine information from the underlying linear formulation to construct a face of the positive semidefinite cone containing the lifted feasible set.
Affine FR can be applied automatically using the linear description of the underlying formulation, and the required affine-hull information is obtained from a single LP followed by linear-algebraic computations. Whether the procedure yields a nontrivial reduction in matrix order, and the magnitude of that reduction, depend on the affine structure of the formulation.

More specifically, this paper makes three contributions. First, we give an explicit affine FR procedure and show that the relaxation it generates satisfies Slater's condition whenever $\aff P=\aff F$. Second, we relate affine FR to analytical facial reduction and, for the Shor relaxation, establish explicit matrix-order comparisons among affine FR, partial FR, and Sieve-SDP. Third, the computational study separately evaluates preprocessing time and reductions in matrix order on MIPLIB mixed-binary instances, and examines SDP solution time and solver reliability on weighted SAT instances.

\paragraph{Notation.} We use $\mathbb{S}^{n}$, $\mathbb{S}_{+}^{n}$ and $\mathbb{S}_{++}^{n}$ to denote the sets of $n \times n$ real symmetric matrices, positive semidefinite matrices, and positive definite matrices, respectively. For two matrices of the same size, the notation $\langle \cdot , \cdot \rangle$ refers to their Frobenius inner product. For a vector $u = (u_1, \ldots, u_n)$, we denote its componentwise absolute value by $|u| = (|u_1|, \ldots, |u_n|)$. The standard basis vectors in $\Rn$ are denoted by $e_{1},\ldots,e_{n}$. When lifting a set from $\Rn$ to $\R^{n+1}$, we occasionally abuse notation and refer to the standard basis vectors in $\R^{n+1}$ as $e_{0},\ldots,e_{n}$. Similarly, the rows and columns of lifted matrices in $\mathbb{S}^{n+1}$ are indexed by $0,1,\ldots,n$. We denote the all-ones and all-zeros vectors of length $k$ by $\ones_{k}$ and $\zeros_{k}$, respectively. The $k \times l$ matrix of all-zeros is denoted by $\textbf{O}_{k \times l}$; when $k = l$, we write $\textbf{O}_{k}$. The $k \times k$ identity matrix is denoted by $I_{k}$. When the dimensions are clear from the context, we omit the subscripts and simply write $\ones,\zeros,\textbf{O}$ and $I$. We use $\Diag(u)$ to denote the diagonal matrix whose diagonal entries are given by the components of $u$. For a subset $S$ of $\R^n$ or $\mathbb{S}^n$, we use $\intr(S)$, $\relint(S)$, and $S^{\perp}$ to denote its interior, relative interior, and the orthogonal complement of the span of $S$, respectively. The convex hull, affine hull, dimension, and cardinality of $S$ are denoted by $\conv(S)$, $\aff(S)$, $\dim(S)$, and $|S|$, respectively. For a closed convex cone $\K$, we write $\K^*$ for its dual cone. Given a linear operator $\mathcal{A}$, we use $\mathcal{A}^*$ to denote its adjoint. The intended meaning of the star notation will be clear from context.

\section{Preliminaries}
\subsection{SDP relaxations for binary optimization problems}
Binary optimization problems aim to find optimal solutions involving binary decision variables. To simplify the presentation, we focus on problems whose feasible region $F$ is the intersection of the binary hypercube and a polyhedron $P$, that is,
\begin{equation}\label{bip}
	F = P \cap \{0,1\}^{n}.
\end{equation}
We assume that the polyhedron $P$ is defined by a linear system $Ax \leq b$, that is, $P = \{ x \in \R^{n} \mid Ax \leq b\}$. Although the comparison in Section~\ref{sec_compare} is restricted to binary programs, the affine FR construction itself uses only the inclusion $F\subseteq P$ and therefore applies without change to mixed-binary feasible sets; see Section~\ref{sec_bip}. Throughout, we assume that $P$ is nonempty and that each binary variable $x_{i}$ is between zero and one for any $x \in P$.

\begin{assumption}\label{ass1}
	$P \neq \emptyset$. If $x \in P$, then $\zeros \leq x \leq \ones$.
\end{assumption}

Our goal is to maximize or minimize a given objective function over the set $F$. In many applications, the objective function is linear or quadratic, in which case SDP relaxations are commonly used, as they can yield strong bounds on the optimal value. SDP is a class of convex optimization problems with linear objective functions and affine constraints over the cone of positive semidefinite matrices. To construct an SDP relaxation for the feasible set $F$, we consider a lifted feasible set
\begin{equation}\label{liftedF}
	F_{1} := \left\{ \begin{bmatrix}1\\x\end{bmatrix}\begin{bmatrix}1\\x\end{bmatrix}^{T} \;\middle|\; x \in F\right\} \subseteq \mathbb{S}_{+}^{n+1}.
\end{equation}
We introduce a matrix variable $Y \in \mathbb{S}^{n+1}$ to relax the nonlinear expression $\left[\begin{smallmatrix}
	1 \\x
\end{smallmatrix}\right]\left[\begin{smallmatrix}
	1 \\x
\end{smallmatrix}\right]^{T}$ in $F_{1}$. Clearly, the matrix variable $Y$ must be positive semidefinite. Additional valid linear equality and inequality constraints on $Y$ can be added to strengthen the relaxation. In general, if $L \subseteq \Sno$ is a nonempty polyhedron such that $F_{1} \subseteq L$, then $F_{1} \subseteq L \cap \Snop$ and the intersection
\[
L \cap \Snop
\]
defines an SDP relaxation for $F_{1}$. By the natural one-to-one embedding of $F$ into $F_{1}$, we also refer to $L \cap \Snop$ as an SDP relaxation of $F \subseteq \mathbb{R}^{n}$.

In practice, many SDP relaxations exist, each offering different tradeoffs between tightness and computational efficiency. One of the most widely used is the \emph{Shor relaxation}, which we describe in detail in \eqref{shorSDP2}. The design of a good SDP relaxation involves balancing its accuracy against the computational effort required to solve it. For a broad discussion on the construction of SDP relaxations for combinatorial optimization problems, see \cite{poljak1995recipe,wolkowicz2000semidefinite}. 

The Shor relaxation is particularly relevant for \emph{quadratically constrained quadratic programs (QCQPs)}, of which the mixed-binary linear programs considered in this paper form a special case. A comprehensive comparison of various SDP relaxations for QCQPs is presented in \cite{bao2011semidefinite}, and the exactness of these relaxations, including the Shor relaxation, is studied in \cite{wang2022quadratically}.

While SDP relaxations can provide strong bounds for binary optimization problems, they come with practical challenges. In particular, generic SDP solvers scale less favorably with problem size than LP solvers and can be computationally expensive. Moreover, due to the specific structure of binary formulations, SDP relaxations may fail Slater's condition, with their feasible regions contained in proper faces of the positive semidefinite cone. A theoretical framework for addressing these issues is discussed in Section \ref{sec_FA}.

\subsection{Shor relaxation} \label{sec:Shor}
The proposed method, affine FR, can be applied to any SDP relaxation of the lifted binary set; see Subsection \ref{sec_affineFRproc}. The face used by affine FR is constructed from an affine set containing the binary feasible points, and therefore does not depend on the particular SDP relaxation chosen. Our theoretical analysis in Section \ref{sec_compare} focuses on the Shor relaxation \cite{shor1987quadratic}; we recall the particular formulation used in this paper below. 
Let $F$ be the binary feasible set defined in \eqref{bip}. Following the indexing convention above, write
\[
Y=(Y_{ij})_{i,j=0}^{n}.
\]
The Shor relaxation of $F$ is given by $L \cap \Snop$, where $L$ is defined as follows:
\linelabel{opA_start}
	\begin{equation}\label{shorSDP2}
		L:= \left\{Y \in \mathbb{S}^{n+1} \, \middle| \, \arrow(Y) = e_{0}, \; \mathcal{A}(Y) \leq \zeros  \right\}.
	\end{equation}
	Here, $e_{0} = (1, 0, \ldots, 0)^T \in \mathbb{R}^{n+1}$ is the first standard basis vector. The lifted operator $\mathcal{A}: \mathbb{S}^{n+1} \rightarrow \R^{m}$ encodes the linear constraints of the original problem and is defined by 
	\begin{equation}\label{shorith}
	    \left(\mathcal{A}(Y)\right)_{i} :=  \left\langle \begin{bmatrix}
		-b_{i} &  \frac{1}{2}a_{i}^{T} \\
		\frac{1}{2}a_{i} & \mathbf{O}
	\end{bmatrix}, Y \right\rangle \text{ for } i = 1,\ldots,m,
	\end{equation}
	\linelabel{opA_end} where $a_{i}^{T}x  \leq  b_{i}$ is the $i$-th constraint in the system $Ax \leq b$. The arrow operator $\arrow : \mathbb{S}^{n+1} \rightarrow \R^{n+1}$ is given by
	\begin{equation}\label{def_arrow}
		\arrow(Y) := \begin{bmatrix}
			Y_{00} \\
			Y_{11} - \frac{1}{2}(Y_{01} + Y_{10}) \\
			\vdots \\
			Y_{nn} - \frac{1}{2}(Y_{0n} + Y_{n0}) \\
		\end{bmatrix} \in \mathbb{R}^{n+1}.
	\end{equation}
	The arrow constraint enforces $Y_{00} = 1$ and $Y_{0i} = Y_{ii}$ for all $i=1, \ldots, n$, serving as a relaxation of the nonlinear binary condition $x_{i} = x_{i}^{2}$ from \eqref{bip}. Accordingly, the matrix variable $Y$ in this relaxation takes the form
	$$Y=\begin{bmatrix} 
		1 & x^{T} \\
		x & X
	\end{bmatrix} \text{ for some } x \in \R^n, X \in \mathbb{S}^{n} \text{ such that } x_{i} = X_{ii} \text{ for } i = 1, \ldots, n.$$

\subsection{Facial reduction}\label{sec_FA}
We say that \emph{Slater's condition} holds for the set $L \cap \Snp$ if $L \cap \intr(\Snp) \neq \emptyset$, i.e., there exists a positive definite feasible solution $Y \in L \cap \Snp$. Slater's condition is a standard regularity condition for conic duality and perturbation analysis. Its failure can contribute to poor numerical behavior, although strict feasibility alone does not guarantee numerical stability.  \linelabel{slater_def}In such cases, we seek an equivalent reformulation of \( L \cap \mathbb{S}_+^n \) that does satisfy Slater's condition. If such a reformulation is found, we say that Slater's condition has been restored.

Facial reduction is a theoretical framework for restoring Slater's condition originally developed by Borwein and Wolkowicz in \cite{borwein1981regularizing,borwein1981facial}. Practical implementations of facial reduction are called \emph{facial reduction algorithms} (FRAs), and various FRAs have been proposed. Although the underlying theory is well established, the computation of suitable exposing vectors remains a central practical difficulty.

\linelabel{def_exp}An essential concept in facial reduction is that of \emph{exposing vectors}, which we now explain:
\begin{enumerate}
    \item If $W \in \Snp$, then the set \begin{equation}\label{expface}
    \{ Y \in \Snp \vs \langle W,Y \rangle = 0\}
\end{equation} is a face of $\Snp$, and $W$ is called an \emph{exposing vector} of this face.\footnote{While the variables in SDP are matrices, the term ``vectors'' remains standard terminology.}

    \item A matrix $W \in \Snp$ is an \emph{exposing vector for} $L \cap \Snp$, if the face exposed by $W$ contains the feasible region, i.e.,\begin{equation}\label{expvec}
	\left\langle W,Y \right\rangle = 0 \quad \forall \;
	Y \in L \cap \Snp.
    \end{equation}
    In particular, if $W$ has the maximum rank among such matrices, it exposes the minimal face of $\Snp$ containing $L \cap \Snp$.
\end{enumerate}

In \cite{borwein1981regularizing,borwein1981facial}, the authors construct exposing vectors using the hyperplane separation theorem. Specifically, given a polyhedron $L$ and a closed convex cone $\K \subseteq \Sn$, the condition $L \cap \relint(\K) = \emptyset$ holds if and only if there exists a hyperplane separating $L$ and $\K$ properly and not containing $\K$, see \cite[Theorem~20.2]{rockafellar1997convex} and \cite[Lemma~3.2]{waki2013facial}. Formally, we have $L \cap \relint(\K) = \emptyset$ if and only if there exists a nonzero matrix $W \in \mathbb{S}^{n}$ such that
\begin{align}
    \label{sepW}
    &\langle W, X \rangle \leq \langle W, Y \rangle && \text{for all } X \in L \text{ and } Y \in \K, \\
    \label{sepWpos}
    &\langle W, Y \rangle > 0 && \text{for some } Y \in \K.
\end{align}
Here, \eqref{sepW} characterizes separation of $L$ and $\K$ by a hyperplane with normal $W$, while \eqref{sepWpos} guarantees that the hyperplane does not contain $\K$. Denote by
\[
\SEP(L \cap \K)
\]
the set of all matrices $W \in \Sn$ that satisfy \eqref{sepW}.\footnote{In this section, we focus on the case where $\K = \Snp$. However, we define $\SEP(L \cap \K)$ for a general closed convex cone $\K$, as this more general form will be used in the analyses in Sections \ref{sec:SEP} and \ref{secPFR}.} We do not include \eqref{sepWpos} in this definition, because $\SEP(L \cap \K)$ will be used below mainly through the homogeneous condition \eqref{sepW}. In the applications considered in this paper, \eqref{sepWpos} is imposed separately when needed; for example, in the positive semidefinite case $\K = \Snp$, it can be enforced by a normalization such as $\operatorname{trace}(W)=1$.

Since $\Snp$ is full-dimensional, every nonzero matrix in $\SEP(L \cap \Snp)$ is an exposing vector for $L \cap \Snp$, and the face exposed by $W$ is a proper face of $\Snp$. However, an exposing vector for $L \cap \Snp$ is not necessarily in $\SEP(L \cap \Snp)$.

\linelabel{Frob_l}
An exposing vector allows us to regularize $L \cap \Snp$. For example, 
if $W = e_1e_1^T$ is an exposing vector (where $e_1$ is the first standard basis vector), then equation~\eqref{expvec} implies that the first row and column of any feasible $Y \in L \cap \Snp$ must be zero. This allows us to reduce the order of the matrix variable: we can replace $Y \in \mathbb{S}_{+}^{n}$ with a smaller matrix variable $R \in \mathbb{S}_{+}^{n-1}$ by removing the first row and column, and reformulate the problem accordingly.

More generally, if $W$ is any nonzero exposing vector for $L \cap \Snp$, we can reformulate $L \cap \Snp$ equivalently as
\begin{equation}\label{K1fr}
    \tilde{L} \cap \Srp \text{ where } \tilde{L} := \left\{R \in \mathbb{S}^{r} \mid VRV^{T} \in L \right\}.
\end{equation}
Here, \( r \) is the dimension of the null space of \( W \), and \( V \in \mathbb{R}^{n \times r} \) is a matrix whose columns form a basis for that null space. We refer to \( V \) as the \emph{facial range vector} associated with the face exposed by \( W \).

Standard facial-reduction theory shows that $L \cap \Snp$ satisfies Slater's condition if and only if $\SEP(L \cap \Snp) = \{\textbf{O}\}$. If a nonzero $W \in \SEP(L \cap \Snp)$ exists, we can use its associated facial range vector $V$ to reduce the problem size via \eqref{K1fr}. The process of identifying such a $W$ and reformulating the problem is called an \emph{FR step}. Repeating this step until no further exposing vectors can be found produces a final reformulation that satisfies Slater's condition.

Assuming that the feasibility problems over $\SEP(L \cap \Snp) \setminus \{\mathbf{O}\}$ can be solved, Slater's condition can always be restored in at most $n$ steps. At each step, choosing a maximum-rank exposing vector maximizes the reduction in the order of the matrix variable among the exposing vectors available at that step. The minimum number of FR steps needed for an FRA to restore Slater's condition is known as the \emph{singularity degree} of $L \cap \Snp$ \cite{sturm2000error}. This quantity is also crucial for deriving error bounds for linear matrix inequalities and plays an important role in conic optimization \cite{drusvyatskiy2015coordinate,im2023revisiting,lin2022optimal,lourencco2021amenable,pataki2018positive,pataki2024exponential,sremac2021error,tanigawa2017singularity}.

Note that the singularity degree of $L \cap \Snp$ is not always one. Even a maximum-rank element of $\SEP(L \cap \Snp)$ may not expose the minimal face containing $L \cap \Snp$, as this set may not contain all exposing vectors. Hence, a single FR step may not restore Slater's condition.

A key challenge in the original FRA proposed by Borwein and Wolkowicz is that finding a nonzero element in $\SEP(L \cap \Snp)$ can be as difficult as solving the SDP itself. More practical approaches use a subset of $\SEP(L \cap \Snp)$, in which finding exposing vectors is easier. In these cases, the resulting reformulation may still fail to satisfy Slater's condition. Nevertheless, such reductions may still reduce the order of the matrix variable and improve numerical behavior in practice. Over the years, various FRAs have been developed using different tractable strategies for identifying exposing vectors. In Section \ref{autoFRA}, we introduce affine FR, which constructs a face from affine information in the underlying binary formulation. 

\subsection{Explicit expression for \texorpdfstring{$\SEP(L \cap \K)$}{SEP}} \label{sec:SEP}
In \cite{borwein1981regularizing,borwein1981facial} and the broader facial reduction literature, it is shown that when $L$ is an affine set, the set $\SEP(L \cap \K)$ equals $L^\perp \cap \K^*$, where $\K \subseteq \Sn$ is a closed convex cone. We next derive an explicit expression for $\SEP(L \cap \K)$ when $L$ is a general polyhedron. We also derive a more specialized expression for the case where $L$ corresponds to the polyhedron in the Shor relaxation, as defined in \eqref{shorSDP2}. 

\begin{lem} \label{lem:SEPpoly}
Let $\K \subseteq \Sn$ be a closed convex cone and $L = \left\{X \in \Sn \mid \mathcal{M}(X) \leq b \right\}$ for some linear operator $\mathcal{M}: \Sn \rightarrow \Rm$ and $b \in \Rm$. Suppose $L \cap \K \neq \emptyset$. Then,
\begin{equation}\label{FRineq}
    \SEP(L \cap \K) = \left\{ \mathcal{M}^{*}(y) \in \K^* \mid y \geq \zeros, \, b^{T}y = 0 \right\}.
\end{equation}

\end{lem}
\begin{proof}
Let $\bar{X} \in L \cap \K$, and recall that $W \in \SEP(L \cap \K)$ if and only if it satisfies \eqref{sepW}. For any $W \in \SEP(L \cap \K)$, we will show that $W \in \K^*$ and $W = \mathcal{M}^{*}(y)$ for some $y \geq \zeros$ with $b^{T}y = 0$.

First, suppose for contradiction that $W \notin \K^*$. Then there exists $\bar{Y} \in \K$ such that $\langle W, \bar{Y} \rangle < 0$. However, \eqref{sepW} implies that
\[
    t \langle W, \bar{Y} \rangle = \langle W, t \bar{Y} \rangle \geq \langle W, \bar{X} \rangle \quad \text{for all } t > 0,
\]
which becomes unbounded below as $t \to \infty$, contradicting the finiteness of $\langle W, \bar{X} \rangle$. Hence, $W \in \K^*$.

Next, since $\mathbf{O}\in \K$, condition \eqref{sepW} with $Y=\mathbf{O}$ gives
\[
   \langle W, X \rangle \leq \langle W,\mathbf{O}\rangle = 0 \quad \text{for all } X \in L.
\]
By LP strong duality, we have 
$$0 \geq \max\left\{ \langle W,X\rangle \mid X \in L \right\} = \min \left\{ b^{T}y \mid \mathcal{M}^{*}(y) = W, \, y \geq \zeros \right\}.$$
Therefore, there exists $y \geq \zeros$ such that $b^{T}y \leq 0$ and $W = \mathcal{M}^{*}(y)$. Moreover, since $\mathcal{M}(\bar{X}) \leq b$, $y \geq \zeros$, $\bar{X} \in \K$, and $W \in \K^*$, 
$$b^{T}y \geq \langle \mathcal{M}(\bar{X}), \, y \rangle = \langle \bar{X} , \mathcal{M}^{*}(y) \rangle  = \langle \bar{X} , W \rangle \geq 0,$$
which implies $b^{T}y = 0$. Hence, $W \in \K^*$ is in the set on the right-hand side of \eqref{FRineq}.

Conversely, suppose $W = \mathcal{M}^{*}(y) \in \K^*$ for some $y \geq \zeros$ with $b^{T}y = 0$. For any $X \in L$ and $Y \in \K$,  
$$ \langle W,X \rangle = \langle \mathcal{M}^{*}(y), X \rangle = \langle \mathcal{M}(X), \,y\rangle \leq b^T y = 0 \leq \langle W,Y \rangle.$$
Therefore, $W \in \SEP(L \cap \K)$.
\end{proof}

When $L$ corresponds to the polyhedron in the Shor relaxation, the expression for $\SEP(L \cap \K)$ can be further specified. 

\begin{lem}\label{shorAUX}
Let $L$ be the polyhedron defined in \eqref{shorSDP2}, and let $\K \subseteq \Sno$ be a closed convex cone. Suppose $L \cap \K \neq \emptyset$. Then,
    \begin{equation}\label{aux}
	\SEP(L \cap \K) = \left\{ \begin{bmatrix}
			-b^{T}y & \frac{1}{2}(y^{T}A - z^{T})  \\
			\frac{1}{2}(A^{T}y - z) & \Diag(z)\\
		\end{bmatrix} \in \K^{*} \, \middle| \, y \in \R^m_+, \ z \in \Rn \right\}.
	\end{equation}

\end{lem}
\begin{proof}
Applying Lemma \ref{lem:SEPpoly} to the polyhedron $L$, where
\[
\mathcal{M} = \begin{bmatrix}
   \mathcal{A} \\ \arrow \\ -\arrow
\end{bmatrix} \ \text{ and the right-hand side is } \ 
\begin{bmatrix}
    \zeros \\ e_0 \\ -e_0
\end{bmatrix},
\]
we obtain that
\begin{align*}
    & \SEP(L \cap \K) \\
    = & \left\{ \mathcal{A}^{*}(y)+\arrow^{*}(z^+)-\arrow^*(z^-) \in \K^{*} \middle| \begin{bmatrix} \zeros \\ e_0 \\ -e_0 \end{bmatrix}^T \begin{bmatrix} y \\ z^+ \\ z^- \end{bmatrix} = 0,  \begin{bmatrix} y \\ z^+ \\ z^- \end{bmatrix} \geq \zeros \right\} \\
    = & \left\{ \mathcal{A}^{*}(y)+\arrow^{*}(\hat{z}) \in \K^{*} \mid e_{0}^{T} \hat{z} = 0, \, \hat{z} \in \mathbb{R}^{n+1}, \, y \in \mathbb{R}_{+}^{m}  \right\}, 
\end{align*}
where we define $\hat{z} = z^+ - z^-$. Writing $\hat{z} = (z_0, z) \in \R \times \R^n$, the constraint $e_{0}^{T}\hat{z} = 0$ implies $z_0 = 0$. Moreover, it is clear that 
\[
\arrow^{*}(\hat{z}) = \begin{bmatrix}
		z_0 & -\frac{1}{2}z^T \\
		-\frac{1}{2}z & \Diag(z)
	\end{bmatrix} \ \text{ and } \ 
   \mathcal{A}^{*}(y) = \begin{bmatrix}
		-b^{T}y & \frac{1}{2}y^{T}A \\
		\frac{1}{2}A^{T}y & \textbf{O}\\
	\end{bmatrix}.
\]
Substituting $z_0 = 0$ and combining the two terms yields the desired expression \eqref{aux}.
\end{proof}

\subsection{Polyhedral theory and linear algebra}
Let $P = \{ x \in \R^{n} \mid Ax \leq b\}$ be a polyhedron. An inequality of the form $c^{T}x \leq \delta$ is said to be \emph{valid} for $P$ if it holds for every $x \in P$.  It is well known that the set of valid inequalities for $P$ can be characterized by nonnegative combinations of the inequalities in $Ax \leq b$.
\begin{lem}[{\cite[Theorem 3.22]{conforti2014integer}}]\label{valid}%
An inequality $c^{T}x \leq \delta$ is \emph{valid} for $P$ if and only if there exists $u \geq \zeros$ such that $u^{T}A = c^{T}$ and $u^{T}b \leq \delta$.
\end{lem}



 
 The $i$-th inequality $a_{i}^{T}x \leq b_{i}$ in the system $Ax \leq b$ is called an \emph{implicit equality} if $a_{i}^{T}x = b_{i}$ holds for every $x \in P$. The \emph{affine hull} of $P$ is the smallest affine set containing $P$. Denoting by $A^{=}x = b^{=}$ the system of all implicit equalities in $Ax \leq b$, the affine hull of $P$ is given by
\begin{equation}\label{affF}
	\aff P = \{ x \in \R^{n} \mid A^{=}x = b^{=}\}.
\end{equation}
The affine-hull description above will be used in the construction of affine FR in the next section.

\section{Affine FR}\label{autoFRA}

In this section, we introduce \emph{affine FR} as a preprocessing method for reducing and, under a sufficient condition, regularizing SDP relaxations for binary optimization problems. The method requires one LP and subsequent linear-algebraic computations, consistent with the \emph{simple and quick} principle for preprocessing advocated by Andersen and Andersen \cite{andersen1995presolving}. Its purpose is to obtain reductions in matrix order at a modest preprocessing cost.

\subsection{The affine FR procedure}\label{sec_affineFRproc}

Affine FR uses exposing-vector language in a slightly different way from classical facial reduction. Classical facial reduction seeks a face containing the whole SDP relaxation $L\cap\Snop$, whereas affine FR first constructs a face containing the lifted binary set $F_1$: a matrix $W\in\Snop$ exposes such a face if $\langle W,Y\rangle=0$ for every $Y\in F_1$. This distinction is important because restricting the SDP relaxation to this face preserves all rank-one matrices corresponding to feasible binary points, while the resulting face restriction may strengthen the chosen relaxation.
Affine FR exploits the fact that any proper affine subset in $\R^n$ containing the feasible region $F \subseteq \{0,1\}^{n}$ defines such a face of the positive semidefinite cone. To see this, suppose that $F \subseteq H$, where $H$ is an affine subspace of dimension $(r < n)$ defined by
\begin{equation}\label{affinel}
H :=\left\{x \in \Rn \mid U^{T} [\begin{smallmatrix} 1\\x \end{smallmatrix}]  = \zeros \right\}, \text{ where } U \in \mathbb{R}^{(n+1) \times (n-r)}.
\end{equation}
Let $W := UU^{T} \in \mathbb{S}_{+}^{n+1}$. Since $U$ has full column rank, we have $\rank W = n-r$. We show that this matrix $W$ satisfies $\langle W,Y \rangle = 0$ for every matrix $Y \in F_{1}$.
Specifically, for any $x \in F$, let $\tilde{x} = [\begin{smallmatrix} 1\\x \end{smallmatrix}]$. It holds that
\begin{equation}\label{wy0}
    \langle W, Y\rangle = \tilde{x}^{T} UU^{T}\tilde{x} =  \norm{U^{T}\tilde{x} }^{2} = 0 \text{ for every } Y = \tilde{x}\tilde{x}^{T} \in F_{1}.
\end{equation}
Let $L \cap \Snop$ be any SDP relaxation for $F$. The condition \eqref{wy0} means that the face
\[
   \mathcal{F}_{W}:=\{Y\in\Snop \mid \langle W,Y\rangle=0\}
\]
contains the lifted binary set $F_1$. Hence we may restrict the relaxation to this face by defining
\[
   \bar L:=L\cap\{Y\mid \langle W,Y\rangle=0\}.
\]
Then $\bar L\cap\Snop$ preserves all lifted feasible binary points and is no weaker than the original relaxation $L\cap\Snop$. The role of affine FR is to work directly in the face $\mathcal{F}_{W}$: if the columns of $V$ form a basis for the null space of $W$, then every matrix in $\mathcal{F}_{W}$ can be written as $Y=VRV^T$ with $R\in\mathbb{S}_{+}^{r+1}$. Thus, replacing $Y$ by $VRV^T$ enforces the face restriction and reduces the order of the matrix variable from $(n+1)$ to $(r+1)$.

To determine an affine set $H$ in \eqref{affinel}, we use the linear system of inequalities that define $F$. Since we assume that $F = P \cap \{0,1\}^{n}$, where $P = \{x \in \Rn \mid Ax \leq b\}$, it follows that $F \subseteq P$, and thus $F \subseteq \aff F \subseteq \aff P$. A practical advantage of using $\aff P$ is that affine FR can be applied without needing the explicit form of the SDP relaxation, since it depends only on the linear system $Ax \leq b$. We describe an efficient subroutine to compute $\aff P$ in Subsection \ref{ahcompute}.

Now we provide a complete description of affine FR. Given $F = P \cap \{0,1\}^{n}$, where $P = \{x \in \Rn \mid Ax \leq b\}$, and an SDP relaxation
\[
   L\cap\Snop,\qquad L:=\{Y\in\Sno\mid \mathcal{M}(Y)\le d\},
\]
where $\mathcal{M}:\Sno\rightarrow\mathbb{R}^{q}$ is a linear operator and $d\in\mathbb{R}^{q}$, the affine FR procedure is as follows.

\begin{enumerate}
    \item \textbf{Compute the affine hull:} Find a matrix $U$ such that $H = \aff P$ is represented in the form given by \eqref{affinel}. 
    \item \textbf{Construct a facial range vector:} Compute a matrix $V$ whose columns form a basis for the null space of $U^T$. This matrix is the facial range vector associated with the positive semidefinite face exposed by $W = UU^T$.
    \item \textbf{Generate the reduced SDP relaxation:} Use the substitution $Y=VRV^T$ to obtain
    \begin{equation}\label{afr_output}
       \tilde L\cap\Srop
        :=
        \left\{R\in\Srop \mid \mathcal{M}(VRV^T)\le d\right\}.
    \end{equation}
   Since $VRV^T\in\mathcal{F}_{W}$ for every $R\in\Srop$, this formulation represents the restriction of the original relaxation to the face exposed by $W$.
\end{enumerate}

We make the following clarifications about the scope and interpretation of affine FR.
\begin{itemize}
    \item \textbf{Applicability.} Affine FR is not tied to any particular SDP relaxation. In the notation above, the chosen relaxation enters only through the linear operator $\mathcal{M}$ and the vector $d$. The construction of $U$, $W$, and $V$ uses the affine information in the original binary formulation and identifies a face of $\Snop$ containing $F_1$; it does not depend on the specific choice of $\mathcal{M}$. Thus, changing from the Shor relaxation to a relaxation with RLT, Lov\'asz--Schrijver, doubly nonnegative, or other valid SDP constraints only changes the operator $\mathcal{M}$ in the reduced constraints $\mathcal{M}(VRV^T)\le d$.
    \item \textbf{Equivalence versus strengthening.} The SDP relaxation generated by affine FR is not necessarily equivalent to the original SDP relaxation. Specifically, it may be strictly stronger if the face restriction exposed by $W$ was not already implied by the original SDP relaxation. Unlike iterative FRAs, affine FR applies this face restriction only once; if the restriction is already implied by the original SDP relaxation, the relaxation generated by affine FR is equivalent to the original SDP relaxation. 
    \item \textbf{Slater restoration.} We say that affine FR \emph{restores Slater's condition} if the relaxation generated by affine FR satisfies Slater's condition. A sufficient condition for affine FR to restore Slater's condition is that the affine set used to construct the face coincides with $\aff F$. In particular, since affine FR uses the computable affine set $H=\aff P$, Slater's condition is restored whenever $\aff P=\aff F$; equivalently, in the notation above, this occurs when $r=\dim(F)$. In this case, the face exposed by $W$ is the minimal face of $\Snop$ containing $\bar{L}\cap\Snop$, and the relaxation generated by affine FR satisfies Slater's condition; see Lemma~\ref{affFisaffP} and \cite{tunccel2001slater}. This is also the viewpoint behind analytical FRAs discussed in Section~\ref{sec:analytical}. When this sufficient condition fails, affine FR may not restore Slater's condition, but it can still reduce the order of the matrix variable and may improve numerical behavior.

\end{itemize}

We provide a concrete example to clarify affine FR next. 
\begin{ex} \label{eg_FR}
	Consider the binary feasible set
	$$F:= P \cap \{0,1\}^{3} = \{ (1,0,0),(0,1,0)\},$$
	where $P = \{ x \in [0,1]^{3} \mid 2x_{1}+x_{2} \leq 2, x_{1}+2x_{2} \leq 2, x_{1}+x_{2} \geq 1, x_{3} \leq 0\}$. The affine hull of $P$ is $$\aff P:=\{ x \in \mathbb{R}^{3} \mid  x_{3} = 0 \}.$$
	The set $F$ and $\aff P$ are depicted in Figures \ref{g1} and \ref{g2}.
	Applying affine FR, we obtain that
	$$U = \begin{bmatrix} 0 \\ 0 \\ 0 \\ 1 \end{bmatrix}, 
        \ W = \begin{bmatrix}
		0\\
		0\\
		0\\
		1
	\end{bmatrix}\begin{bmatrix}
		0\\
		0\\
		0\\
		1
	\end{bmatrix}^{T},  \text{ and } V = \begin{bmatrix}
		1 & 0& 0\\
		0 & 1 & 0\\
		0 & 0 & 1\\
		0 & 0 & 0\\
	\end{bmatrix}.$$
	Thus, we can reduce the order of the matrix variable in an SDP relaxation for $F$ from $4$ to $3$.
	\begin{figure}[H]
		\centering
		\begin{minipage}[t]{0.45\textwidth}
			\centering
			\begin{tikzpicture}[scale=2.5, transform shape, rotate around z = 0]
				
				\draw[color=red,fill=red!100] (1,0,0) circle [radius=.2mm];
				\draw[color=red,fill=red!100] (0,0,1) circle [radius=.2mm];
				
				\coordinate (A) at (0, 0, 0);
				\coordinate (B) at (1.3, 0, 0); 
				\coordinate (C) at (0, 1, 0);
				\coordinate (D) at (0, 0, 1.3);
				\draw[dashed,->,opacity=.5] (A) -- (B);
				\draw[dashed,->,opacity=.5] (A) -- (C);
				\draw[dashed,->,opacity=.5] (A) -- (D);
				
				\draw[fill=gray!20] (1,0,0) -- (0,0,1) -- (2/3,0,2/3) -- (1,0,0);
			\end{tikzpicture}
			\caption{The red dots constitute $F$, and the grey triangular area represents $P$.}\label{g1}
		\end{minipage}\hfill
		\begin{minipage}[t]{0.45\textwidth}
			\centering
			\begin{tikzpicture}[scale=2.5, transform shape, rotate around z = 0]
				
				\draw[color=red,fill=red!100] (1,0,0) circle [radius=.2mm];
				\draw[color=red,fill=red!100] (0,0,1) circle [radius=.2mm];
				
				\coordinate (A) at (0, 0, 0);
				\coordinate (B) at (1.3, 0, 0); 
				\coordinate (C) at (0, 1, 0);
				\coordinate (D) at (0, 0, 1.3);
				\draw[dashed,->,opacity=.5] (A) -- (B);
				\draw[dashed,->,opacity=.5] (A) -- (C);
				\draw[dashed,->,opacity=.5] (A) -- (D);
				
				\fill[fill=gray!20] (1,0,0) -- (1,0,1) -- (0,0,1) -- (0,0,0) -- (1,0,0);
			\end{tikzpicture}
			\caption{The affine hull of $P$ is represented by the (unbounded) gray area.}\label{g2}
		\end{minipage}
	\end{figure}
	
\end{ex}

In the previous example, $\aff P\neq\aff F$, so the sufficient condition for Slater restoration in Lemma~\ref{affFisaffP} does not apply. The continuation below instead chooses $H=\aff F$, which guarantees Slater restoration.

\begin{ex}
The affine hull of $F:= \{ (1,0,0),(0,1,0)\} \subseteq \R^3$ is given by $$\aff F=\{ x \in \mathbb{R}^{3} \mid x_{1}+x_{2} = 1, x_{3} = 0 \}.$$
This is a one-dimensional affine subspace, as illustrated in Figures \ref{g3} and \ref{g4}. From this affine hull, we obtain that
\[
U = \begin{bmatrix} -1 & 0 \\ 1 & 0 \\ 1 & 0 \\ 0 & 1 \end{bmatrix}, \ 
W = UU^T = \begin{bmatrix}
1 & -1 & -1 & 0\\
-1 & 1 & 1 & 0\\
-1 & 1 & 1 & 0\\
0 & 0 & 0 & 1\\
\end{bmatrix}, \text{ and } 
V = \begin{bmatrix}
	1 & 1\\
	1 & 0 \\
	0 & 1 \\
	0 & 0\\
\end{bmatrix}.
\]
Thus, we can reduce the order of the matrix variable in an SDP relaxation from $4$ to $2$, and Slater's condition is restored.

\begin{figure}[H]
	\centering
	\begin{minipage}[t]{0.45\textwidth}
		\centering
			\begin{tikzpicture}[scale=2.5, transform shape, rotate around z = 0]
			
			\draw[color=red,fill=red!100] (1,0,0) circle [radius=.2mm];
			\draw[color=red,fill=red!100] (0,0,1) circle [radius=.2mm];
			
			\coordinate (A) at (0, 0, 0);
			\coordinate (B) at (1.3, 0, 0); 
			\coordinate (C) at (0, 1, 0);
			\coordinate (D) at (0, 0, 1.3);
			\draw[dashed,->,opacity=.5] (A) -- (B);
			\draw[dashed,->,opacity=.5] (A) -- (C);
			\draw[dashed,->,opacity=.5] (A) -- (D);
			
		\end{tikzpicture}
		\caption{The red dots constitute $F$.}\label{g3}
	\end{minipage}\hfill
	\begin{minipage}[t]{0.45\textwidth}
		\centering
	\begin{tikzpicture}[scale=2.5, transform shape, rotate around z = 0]
	
	\draw[color=red,fill=red!100] (1,0,0) circle [radius=.2mm];
	\draw[color=red,fill=red!100] (0,0,1) circle [radius=.2mm];
	
	\coordinate (A) at (0, 0, 0);
	\coordinate (B) at (1.3, 0, 0); 
	\coordinate (C) at (0, 1, 0);
	\coordinate (D) at (0, 0, 1.3);
	\draw[dashed,->,opacity=.5] (A) -- (B);
	\draw[dashed,->,opacity=.5] (A) -- (C);
	\draw[dashed,->,opacity=.5] (A) -- (D);
	
	\draw[gray] (1.1,0,-.1) -- (-.1,0,1.1);
\end{tikzpicture}
		\caption{The affine hull of $F$ is represented by the line connecting the two red points.}\label{g4}
	\end{minipage}
\end{figure}

\end{ex}

\subsection{Computing the affine hull}\label{ahcompute}
In this subsection, we discuss implementation details concerning the computation of the affine hull of a polyhedron $P$. Specifically, to determine $\aff P$, it is necessary to identify all implicit equalities in the system $Ax \leq b$ defining $P$, as described in \eqref{affF}. In \cite{fukuda2016lecture}, Fukuda proposed a method that solves up to $m$ LPs associated with $P$ to detect all implicit equalities. However, this approach can be computationally expensive, particularly when the number of inequalities $m$ is large. Here, we present a more efficient alternative based on a single LP formulation.

Assuming $P \neq \emptyset$, LP strong duality implies that the following primal and dual LPs share the same optimal value:
\begin{equation*}
    \max \{ 0 \mid A x \leq b\} = \min \{ b^{T}y \mid A^T y = \zeros, y \geq \zeros\}.
\end{equation*}
Note that any feasible solution $x$ to the primal is optimal. Let $y^{*}$ be any optimal solution to the dual. By the complementary slackness condition, it holds that  $(b-Ax)^{T}y^{*} = 0$. Let $I := \{ i \mid y_{i}^{*} > 0\}$ be the index set of positive entries in $y^{*}$. The primal constraints in $Ax \leq b$ associated with $I$ are always active, that is, for any $x \in P$, we have $a_{i}^{T}x = b_{i}$ for each $i \in I$. By definition, these constraints are implicit equalities of the system. 

If the dual optimal solution $y^{*}$ has the maximum number of nonzeros, then the corresponding set $I$ captures all implicit equalities by \linelabel{thm_scsl}\label{thm_scs} the Goldman-Tucker Theorem; see \cite[Corollary 2A]{goldman1956theory}. In this case, the affine hull of $P$ is given by
\begin{equation*}
    \aff P = \{x \in \Rn \mid a_{i}^{T}x = b_{i} \text{ for every } i \in I\}.
\end{equation*}
There are several ways to compute such an optimal dual solution $y^{*}$ with the maximum number of nonzeros. One approach is to use interior point methods, as shown by G\"{u}ler and Ye \cite{guler1993convergence}. Alternatively, $y^*$ can be obtained by solving a single LP. We recall the following standard split-variable formulation, based on \cite[Exercise 3.27]{bertsimas1997introduction} and \cite{mehdiloo2021finding}.

Consider the following LP:
\begin{equation}\label{splitFR}
	\begin{array}{rrll}
		\max & \multicolumn{1}{l}{\ones^{T}u} \\
		\text{subject to} & (u+v)^{T}\begin{bmatrix}
			A & b
		\end{bmatrix} &=& \zeros^{T} \\
		& u,v &\geq& \zeros\\
		& u &\leq& \ones.
	\end{array}
\end{equation}
To see why \eqref{splitFR} has the desired effect, let
\[
   K:=\{y\in\mathbb{R}^{m}_{+}\mid A^{T}y=\zeros,\ b^{T}y=0\}.
\]
The vectors in $K$ are precisely the dual optimal solutions. Let $k:=\max\{|\operatorname{supp}(y)|\mid y\in K\}$. For any feasible $(u,v)$ in \eqref{splitFR}, the vector $y=u+v$ belongs to $K$, and therefore
\[
   \ones^{T}u\leq |\operatorname{supp}(u)|\leq |\operatorname{supp}(y)|\leq k.
\]
Conversely, choose $y\in K$ with $|\operatorname{supp}(y)|=k$. Since $K$ is a cone, we can scale $y$ so that all of its positive components are at least one. Setting $u_i=\min\{y_i,1\}$ and $v=y-u$ gives a feasible solution of \eqref{splitFR} with objective value $k$. Hence \eqref{splitFR} identifies the maximum possible support of a dual optimal solution. If $(u^{*},v^{*})$ solves \eqref{splitFR}, then $y^{*}=u^{*}+v^{*}$ is a dual optimal solution with maximum support, and $I=\{i\mid y_i^{*}>0\}$ gives all implicit equalities.

Thus, we can compute the affine hull of a polyhedron with $n$ variables and $m$ inequalities by solving the LP in \eqref{splitFR}, which has $2m$ variables and $n+1$ equality constraints. As the variable $u$ is bounded above, the bounded-variable simplex method can be applied to solve \eqref{splitFR}. Although this formulation doubles the number of variables, it remains a linear program. Section~\ref{sec_bip} reports its observed computational cost in the benchmark experiments.

\begin{remark}
	The costs of solving \eqref{splitFR} can be reduced as follows. Given a feasible solution $x \in P$, we can check whether each inequality $a_{i}^{T}x \leq b_{i}$ holds strictly. If so, it cannot be an implicit equality, and there is no need to include the variables $u_{i}$ and $v_{i}$ in \eqref{splitFR}. In many applications, it is relatively easy to generate feasible solutions, and they can be used to reduce the size of \eqref{splitFR} substantially.
\end{remark}

\section{A Theoretical Comparison}\label{sec_compare}
In this section, we present a theoretical comparison between affine FR and several existing approaches in the literature. Facial reduction involves identifying an exposing vector, as defined in \eqref{expvec}. A central challenge in FRAs lies in efficiently computing such exposing vectors in practice. 
Affine FR itself is independent of the particular SDP relaxation, as discussed in Subsection~\ref{sec_affineFRproc}. In this section, however, we specialize to the Shor relaxation because it provides a canonical common setting in which the exposing vectors and resulting reductions in matrix order of different FRAs can be compared explicitly. Thus, the results below should be read as a theoretical comparison on the Shor relaxation, not as a restriction on the applicability of affine FR.

Existing FRAs for SDP can generally be categorized into two main types:
\begin{enumerate}
	\item \textbf{Analytical methods} use problem-specific formulas to construct exposing vectors, typically from a closed-form description of the affine hull of the feasible set. Their computational cost is negligible once such a description is available, but deriving it requires additional structural information.
	
	\item \textbf{General-purpose methods}, including partial FR \cite{permenter2018partial} and Sieve-SDP \cite{zhu2019sieve}, search for exposing vectors using the chosen SDP representation. Their reductions may therefore depend on the linear description of the relaxation and, for partial FR, on the auxiliary cone used in the search.
\end{enumerate}

The key features and characteristics of these methods are summarized in Table \ref{table_comparemethod}.

\begin{table}[htbp]
	\centering
	\begin{tabular}{|c|c|c|c|c|} \hline
		Method& Scope & Information used & Computation & \shortstack{SDP\\dependent}\\ \hline
		Analytical& Problem-specific & $\aff F$ & Analytical derivation & No\\ \hline
		Sieve-SDP& General SDP & Constraint matrices & Matrix tests & Yes\\ \hline
		Partial FR& General SDP & $\SEP(L\cap\K)$ & Conic problem & Yes\\ \hline
		Affine FR& Binary/mixed-binary & $\aff P$ & One LP & No\\ \hline
	\end{tabular}
	\caption{\enspace Structural comparison of FRAs considered in this paper}
	\label{table_comparemethod}
\end{table}

Each method in Table \ref{table_comparemethod} has advantages in different settings. Analytical methods are particularly attractive when problem-specific structure yields an explicit affine-hull description. Sieve-SDP and partial FR are designed for general SDP formulations, whereas affine FR is tailored to SDP relaxations of binary and mixed-binary optimization problems. These methods should therefore be viewed as complementary, with their applicability and effectiveness depending on the problem class and its formulation.

The following theorem compares the matrix orders obtained by the four methods on the Shor relaxation.

\begin{thm} \label{thm_main}
Let $F:= P \cap \{0,1\}^{n}$, where $P = \{ x\in \Rn \vs Ax \leq b \}$, as defined in \eqref{bip}. We consider four FRAs applied to the Shor relaxation for $F$, as summarized in Table \ref{tab:fra}. Let $r_A$, $r_P$, $r_P^+$, and $r_S$ denote the orders of the matrix variables obtained after applying affine FR, partial FR with $\K=\Ds$, partial FR with $\K=\DDs$, and Sieve-SDP, respectively, to the Shor relaxation for $F$.
Then, 

	$$r_{A} \leq r_{P}^{+} \leq r_{P} \leq r_{S} = n+1.$$
\end{thm}

\begin{table}[htbp]
    \centering

\begin{tabular}{|c|c|}\hline
    Reduced Matrix Order     & Method \\ \hline
    $r_{A}$         & Affine FR (Section \ref{autoFRA}) \\ \hline
   $r_{P}$         & Partial FR with cone $\Ds$ (Section \ref{pfrd}) \\ \hline
   $r_{P}^{+}$     & Partial FR with cone $\DDs$ (Section \ref{sec_ddm}) \\ \hline
    $r_{S}$         & Sieve-SDP (Section \ref{sec_sieve}) \\ \hline
\end{tabular}
    \caption{Notation for the order of the matrix variable after applying FRAs to the Shor relaxation}
    \label{tab:fra}

\end{table}

The proof of Theorem \ref{thm_main} is presented over the remainder of this section, through Corollary \ref{lem_rarp2} and Lemma \ref{lem_sieve}.

\subsection{Analytical methods} \label{sec:analytical}
Classical analytical FRAs can be viewed as choosing $H=\aff F$, whereas affine FR uses the computable affine set $H=\aff P$. Analytical methods derive $\aff F$ in closed form from problem-specific structure, while affine FR computes $\aff P$ from the linear formulation. In \cite{zhao1998semidefinite}, Zhao et al. applied this analytical approach to construct SDP relaxations for the quadratic assignment problem that satisfy Slater's condition. Their method relies on the fact that the set of permutation matrices has a compact convex hull representation:
\begin{equation}\label{AP_QAP}
	\conv F = \{X \in [0,1]^{n \times n} \mid X\textbf{1}=\textbf{1}, X^{T}\textbf{1} = \textbf{1}, X \geq \textbf{O}\}.
\end{equation}
As $X \geq \textbf{O}$ are not implicit equalities, the affine hull of $F$ is fully determined by the equality constraints in \eqref{AP_QAP}. This structure enables an analytical derivation of $\aff F$, facilitating a reduction that restores Slater's condition.

When such a closed-form expression for $\aff F$ can be derived, the analytical approach requires little additional computation. However, its applicability is restricted to problem classes for which the required affine-hull description is available. Most studies employing the analytical approach focus on variants of the quadratic assignment problem, including the graph partitioning problem \cite{lisser2003graph,wolkowicz1999semidefinite}, the min-cut problem \cite{li2021strictly}, and the vertex separator problem \cite{rendl2018min}. Similar techniques have also been used in the quadratic cycle covering problem \cite{de2021sdp} and the quadratic shortest path problem \cite{hu2020solving} to restore Slater's condition for their SDP relaxations. 
In general, however, computing $\aff F$ explicitly is often computationally challenging; in fact, it is NP-hard in general.\footnote{For example, given a set of positive integers $a_{1},\ldots,a_{n}$, the partition problem asks whether there exists a subset $T \subseteq \{1,\ldots,n\}$ such that $\sum_{i\in T}a_{i} = \sum_{i\notin T}a_{i}$. The partition problem is known to be NP-hard; see \cite[page 97]{karp2010reducibility}. Define the binary set $F = \{ x \in \{0,1\}^{n} \vs \sum_{i=1}^{n} a_{i} x_{i} = \sum_{i=1}^{n} a_{i} (1-x_{i}) \}$. The answer to the partition problem is NO if and only if $F = \emptyset$. Determining whether $F = \emptyset$ is equivalent to checking whether the affine hull of $F$ is empty.} This limitation is a key motivation for developing the affine FR approach introduced in this work. 

When applied to the structured binary optimization problems mentioned above, affine FR restores Slater's condition. This outcome follows directly from the fact that $\aff F = \aff P$ holds for all these problems, together with Lemma \ref{affFisaffP} below. Therefore, in situations where it is unclear whether a binary optimization problem admits a tractable analytical formula for its affine hull, affine FR provides a computable alternative based on $\aff P$, subject to the sufficient condition in Lemma \ref{affFisaffP}.

\begin{lem}\label{affFisaffP}
Let $F = P \cap \{0,1\}^{n}$ for some polyhedron $P$. If $\aff F = \aff P$, then for any SDP relaxation $L\cap\Snop$ for $F$, the SDP relaxation $\tilde L\cap\Srop$ generated by affine FR in \eqref{afr_output} satisfies Slater's condition.
\end{lem}
\begin{proof}
Let $L \cap \Snop$ be an arbitrary SDP relaxation for $F$, and let $U$, $W=UU^T$, and $V$ be the matrices generated by affine FR from $H=\aff P$. Since $\aff F = \aff P$, we have $r:=\dim F = \dim P$.

Choose $r+1$ affinely independent points $v_{1},\ldots,v_{r+1} \in F$, and set $\tilde v_i:=\begin{bmatrix}1\\ v_i\end{bmatrix}$. For each $i$, the matrix $\tilde v_i\tilde v_i^T$ belongs to $F_1$, and hence to $L\cap\Snop$ by the definition of an SDP relaxation. Since $L$ and $\Snop$ are convex, the matrix
\[
    \hat{Y} := \frac{1}{r+1}\sum_{i=1}^{r+1} \tilde v_i\tilde v_i^T
\]
also belongs to $L\cap\Snop$. Moreover, the vectors $\tilde v_1,\ldots,\tilde v_{r+1}$ are linearly independent, so $\rank \hat{Y}=r+1$.

Because $H=\aff P=\aff F$, we have $U^T\tilde v_i=\zeros$ for all $i$. Therefore $\langle W,\tilde v_i\tilde v_i^T\rangle=0$ for all $i$, and hence $\hat Y$ lies in the face exposed by $W$. The null space of $W$ has dimension $r+1$, and the columns of $V\in\mathbb{R}^{(n+1)\times(r+1)}$ form a basis for this null space. Thus there exists $\hat R\in\Srop$ such that $\hat Y=V\hat R V^T$. Since $\hat Y\in L$, this matrix $\hat R$ is feasible for the reduced relaxation $\tilde L\cap\Srop$ generated in \eqref{afr_output}.

Finally, since $V$ has full column rank, $\rank(V\hat R V^T)=\rank(\hat R)$. Hence $\rank(\hat R)=\rank(\hat Y)=r+1$. As $\hat R\in\Srop$, this implies $\hat R\succ0$. Therefore $\tilde L\cap\Srop$ contains a positive definite feasible point and satisfies Slater's condition.
\end{proof}

\subsection{Partial FR}\label{secPFR}
Permenter and Parrilo \cite{permenter2018partial} proposed a specialized FRA, which we refer to as \emph{partial FR}. When applied to $L \cap \Snp$, partial FR aims to find, at each FR step, a \emph{maximum-rank} matrix $W \in \SEP(L \cap \K)$, where $\K$ is a closed convex cone containing $\Snp$. Since $\K \supseteq \mathbb{S}_{+}^{n}$, it follows that $\K^{*} \subseteq \mathbb{S}_{+}^{n}$, and thus, by Lemma \ref{lem:SEPpoly}, $\SEP(L \cap \K) \subseteq \SEP(L \cap \Snp)$. Therefore, any nonzero matrix $W$ identified in such an FR step serves as an exposing vector for $L \cap \Snp$. Following the general facial reduction framework described in Section \ref{sec_FA}, partial FR iteratively applies these steps until no further reduction is possible. 

However, because $\SEP(L \cap \K)$ may not contain all exposing vectors for $L \cap \Snp$, partial FR may terminate early, failing to restore Slater's condition. On the other hand, partial FR benefits from low computational cost when $\K$ is chosen appropriately. Following the suggestions in \cite{permenter2018partial}, we consider two specific choices for $\K$:

\begin{enumerate}
	\item The cone $\Ds$ is the set of symmetric matrices with nonnegative diagonal entries: 
    $$\Ds:=\{ X \in \Sn \mid X_{ii} \geq 0 \;\forall\; i \}.$$ 
    Its dual cone $\D$ is the set of diagonal matrices with nonnegative diagonal entries:
	$$\D= \left\{ W \in \Sn \mid  W_{ii} \geq 0, \, W_{ij} = 0 \;\forall\; i \neq j \right\}.$$ 
	\item The cone $\DDs$ is defined as
    \begin{equation}\label{def_dd1}
    \DDs:= \left\{ X \in \Sn \;\middle|\; X_{ii} \geq 0, \, X_{ii} + X_{jj} \pm 2X_{ij} \geq 0 \;\forall\; i,j\right\}.        
    \end{equation}
    Its dual cone $\DD$ is the set of symmetric diagonally dominant matrices:
\begin{equation}\label{def_dd}
    \DD= \left\{ W \in \Sn \;\middle|\; W_{ii} \geq \sum_{j\neq i} |W_{ij}| \;\forall\; i \right\}.
\end{equation}
\end{enumerate}
Since the dimensions of the cones above are always evident from the context, we omit their dimensional superscripts to simplify the notation.

We discuss these two choices in detail in the following subsections. Throughout this subsection, we adopt the same notation as in Section~\ref{sec:Shor}, unless stated otherwise. Specifically, let $F=P\cap\{0,1\}^{n}$, where $P=\{x\in\R^{n}\mid Ax\leq b\}$, and let $L\cap\Snop$ denote the Shor relaxation for $F$, where $L$ is defined in \eqref{shorSDP2}. This relaxation contains the arrow equations and the lifted original linear inequalities, but no equality-generated RLT equations. We assume that Assumption~\ref{ass1} holds throughout this subsection.

This formulation-specific scope is important. Because the lifted constraints $\mathcal{A}(Y)\leq\zeros$ are retained as inequalities, $L$ is a polyhedron in the matrix space rather than an affine subspace. Introducing nonnegative slack variables gives the natural affine conic formulation over the product cone $\Snop\times\mathbb{R}_{+}^{m}$, rather than an affine slice in the original positive semidefinite matrix space considered in \cite{huXu2026face}. Consequently, the results below are not direct special cases of that analysis.

\begin{remark}[Representation dependence of partial FR]
\label{rem:partial_fr_representation}
Partial FR with a restricted cone such as $\Ds$ or $\DDs$ is representation dependent. Indeed, once an exposed face is described by a facial range vector $V$ through the substitution $Y=VRV^T$, the auxiliary cone $\K$ is imposed on the coordinate matrix $R$. Choosing a different facial range vector for the same face, or equivalently applying an invertible congruence transformation to $R$, can change the cones $\Ds$ and $\DDs$, and hence the exposing vectors detected by partial FR. In the analysis below, we do not choose $V$ arbitrarily. The exposed faces have a natural interpretation inherited from the lifting used in the Shor relaxation: variable-indexed rows and columns are retained, set to zero, or identified with the zeroth row and column according to the variables fixed by $P$. This interpretation determines the facial range vector used in our analysis, and each subsequent partial FR step is taken with respect to the corresponding facial range vector and its induced row/column indexing.
\end{remark}

\subsubsection{Non-negative diagonal matrices}\label{pfrd}

We first simplify the expression for $\SEP(L \cap \K)$ when $\K = \Ds$.
\begin{lem}\label{PFRDlem1}
It holds that
    \begin{equation}\label{PFRDclaim2}
              \SEP(L \cap \Ds) = \left\{\begin{bmatrix}
               0 & \textbf{0}^{T} \\
        \textbf{0} & \Diag\left(A^{T} y\right)
        \end{bmatrix} \middle| A^{T}y \geq \zeros, \, b^{T}y = 0 , \, y \in \R^m_+ \right\}.
        \end{equation}
\end{lem}
\begin{proof}
    By Lemma \ref{shorAUX}, we have
    \begin{equation*}
		\SEP(L \cap \Ds) = \left\{ \begin{bmatrix}
			-b^{T}y & \frac{1}{2}(y^{T}A - z^{T})  \\
			\frac{1}{2}(A^{T}y - z) & \Diag(z)\\
		\end{bmatrix} \in \D \middle| \, y \in \R^m_+, z \in \R^n \right\}.
    \end{equation*}
    Since $\D$ consists of diagonal matrices with nonnegative diagonal entries, the above expression simplifies to
    \begin{equation} \label{PFRDlem1eq1}
		\SEP(L \cap \Ds) = \left\{ \begin{bmatrix}
			-b^{T}y & \zeros^T  \\
			\zeros & \Diag(A^T y)\\
		\end{bmatrix} \, \middle| \, b^T y \leq 0, \, A^Ty \geq \zeros, \, y \in \R^m_+\right\}.
    \end{equation}

    By Assumption \ref{ass1}, the polyhedron $P$ is nonempty and $x \geq \zeros$ for every $x \in P$. Fix any $x^* \in P$. For any $y \in \R^m_+$ satisfying $b^T y \leq 0$ and $A^Ty \geq \zeros$, we have that 
    \begin{equation} \label{pfr_eq1}
    0 \leq y^{T}Ax^{*} \leq y^{T}b \leq 0.
    \end{equation}
    As a result, \eqref{PFRDlem1eq1} is equivalent to \eqref{PFRDclaim2}, which completes the proof.
        
\end{proof}

In the next lemma, we present an explicit form of the exposing vectors generated by partial FR with $\K = \Ds$. We also show that partial FR with $\K=\Ds$ finds no further exposing vector after at most one FR step. The key observation underlying the proof is that, after one FR step, the new matrix variable $R$ is obtained from the original matrix variable $Y$ by removing the rows and columns corresponding to the variables that are fixed at zero in $P$.

\begin{lem}\label{PFRD}
For each $i \in \{1, \ldots,n\}$, define
\[
     W^{i} := \begin{bmatrix} 0 & \zeros^{T} \\ 
                \zeros & e_{i}e_{i}^{T}    
                \end{bmatrix} \in \Sno,
\]
where $e_i \in \R^n$ is the i-th standard basis vector. Let 
\[
I_+:= \{ i \in \{1,\ldots,n\} \mid x_{i} = 0 \text{ for all } x \in P\},
\]
and define $W^* = \sum_{i \in I_+} W^i$. 
Then the following statements hold:
\begin{enumerate}
    \item $W^{*} \in \SEP(L \cap \Ds)$.
    \item $W^{*}$ has the maximum rank among all matrices in $\SEP(L \cap \Ds)$.
    \item When applied to $L \cap \Snop$, partial FR with $\K=\Ds$ finds no further exposing vector after at most one FR step.
\end{enumerate}
    
\end{lem}
\begin{proof}
\begin{enumerate}
    \item If $i \in I_+$, then $x_{i} = 0$ for every $x\in P$. LP strong duality implies that
	\begin{equation}\label{pfr_eq2}
		\min \{ b^{T}y \mid A^{T}y = e_{i} ,y \geq \zeros \} =	\max\{ x_{i} \mid Ax \leq b\} =	0.
	\end{equation}
    Therefore, there exists $y \in \R^m_+$ such that $A^Ty = e_i$ and $b^Ty = 0$, which implies $W^i \in \SEP(L \cap \Ds)$ by Lemma \ref{PFRDlem1}. As $\SEP(L \cap \Ds)$ is a polyhedral cone, it follows $W^* = \sum_{i \in I_+} W^i \in \SEP(L \cap \Ds)$.

 \item If $i \notin I_+$, then $x_{i}^*>0$ for some $x^* \in P$. In this case, any $y \in \R^m_+$ satisfying $A^T y \geq \zeros$ and $b^T y = 0$ must also satisfy \eqref{pfr_eq1}, and therefore $(A^{T}y)_{i} = 0$. Consequently, every matrix $W \in \SEP(L \cap \Ds)$ must have zero entries in the $i$-th row and column whenever $i \notin I_+$ (with indices starting from 0). Therefore, the only positions where nonzero entries can appear in $W$ are those indexed by $i \in I_+$. Since $W^*$ is constructed to be nonzero exactly in these positions, it achieves the maximum possible rank among all matrices in $\SEP(L \cap \Ds)$.
 
\item If $\SEP(L \cap \Ds) = \{\mathbf{O}\}$, then partial FR with $\K = \Ds$ terminates immediately. Otherwise, the algorithm generates a maximum-rank exposing vector in $\SEP(L \cap \Ds)$. By the previous item, any maximum-rank element of $\SEP(L\cap\Ds)$ has positive diagonal entries exactly in the positions indexed by $I_+$, and hence has the same null space as $W^*$. Therefore, it is enough to describe the reduction using $W^*$. Let $V \in \R^{(n+1)\times(n+1-|I_+|)}$ be the associated facial range vector, obtained by removing the $i$-th column of the identity matrix $I_{n+1}$ for each $i \in I_+$ (with column indices starting from 0). 
Using the reformulation in \eqref{K1fr}, the matrix variable $R$ in the reduced problem is obtained from $Y$ by removing the rows and columns corresponding to indices $i \in I_+$. Let $\tilde{P}$ be the polyhedron obtained by eliminating from $P$ all variables $x_i$ with $i \in I_+$, and define $\tilde{F}:= \tilde{P} \cap \{0,1\}^{n-|I_+|}$. Then the reformulation $\tilde{L} \cap \mathbb{S}^{n+1-|I_+|}_+$ coincides with the Shor relaxation for $\tilde{F}$. Since $\tilde{P}$ contains no variables fixed at zero, we have $\SEP(\tilde{L} \cap \Ds) = \{\mathbf{O}\}$, and partial FR with $\K=\Ds$ finds no further exposing vector after this step.
    \end{enumerate}
\end{proof}

For an index $i \notin I_+$, it is nevertheless possible that $x_i>0$ for some $x\in P$, while $x_i=0$ for all $x\in F$. Consequently, this approach does not detect all variables that are fixed at zero in the binary feasible set $F$. To compute the maximum-rank exposing vector $W^{*}$ in Lemma \ref{PFRD}, one may apply the technique described in Subsection \ref{ahcompute}. An alternative method is proposed in \cite{permenter2018partial}.



\begin{cor} \label{lem_rarp}
It holds that $r_{A} \leq r_{P}$.
\end{cor}
\begin{proof}

By definition of the index set $I_+$, it holds that $P \subseteq \{x \in \R^n \mid x_i = 0 \text{ for } i \in I_+\}$. Therefore, $\dim(P) \leq \dim(\{x \in \R^n \mid x_i = 0 \text{ for } i \in I_+\}) = n-|I_+|$. By construction, $r_A = \dim(P) + 1$ and $r_P = n+1 - |I_+|$. Therefore, $r_A \leq r_P$, as claimed.
\end{proof}

We illustrate the ideas presented in the preceding results with the following examples.
\begin{ex} \label{eg:pFRD}
	Let $P$ be the polyhedron defined by the linear system $Ax \leq b$, where
    \[A = \begin{bmatrix}
			1 & 1 \\
			-1 & -1\\
			1 & 0 \\
			-1 & 0 \\
			0 & -1 \\
	\end{bmatrix} \text{ and }  b = \begin{bmatrix}
			1\\
			-1\\
			0\\
			0\\
			0
	\end{bmatrix}.
    \]
    It is clear that $P = \{(x_1, x_2) \in \R^2 \mid x_1 + x_2 = 1, \, x_1 = 0, \, x_2 \geq 0\} =  \{(0,1)\}$. Let $F:= \{0,1\}^{2} \cap P = \{(0,1)\}$. The corresponding Shor relaxation is $L \cap \mathbb{S}_{+}^{3}$, where
    \begin{align*}
    L &= \left\{Y = \begin{bmatrix}
        Y_{00} & Y_{01} & Y_{02} \\
        Y_{10} & Y_{11} & Y_{12} \\
        Y_{20} & Y_{21} & Y_{22}
    \end{bmatrix} \in \mathbb{S}^3 \, \middle| \, 
    \begin{array}{c}
         Y_{00} = Y_{01}+Y_{02}, \, Y_{01} = 0, \, Y_{02} \geq 0,   \\
         Y_{00} = 1, \, Y_{11} = Y_{01}, \, Y_{22} = Y_{02}
    \end{array}
    \right\} \\
    &= \left\{\begin{bmatrix}
        1 & 0 & 1 \\
        0 & 0 & Y_{12} \\
        1 & Y_{12} & 1
    \end{bmatrix} \, \middle| \,  Y_{12} \in \R
    \right\}.
    \end{align*}

    Now apply partial FR with $\Ds$ to $L \cap \mathbb{S}_{+}^{3}$. Lemma \ref{PFRDlem1} implies that
    \begin{align*}
    \SEP(L \cap \Ds) 
    & = \left\{ \begin{bmatrix}
        0 & 0 & 0 \\
	0 & y_1-y_2+y_3-y_4 & 0 \\
	0 & 0 & y_1-y_2-y_5
	\end{bmatrix} \middle| \begin{array}{l}
	y_1-y_2+y_3-y_4 \geq 0, \\
        y_1-y_2-y_5 \geq 0, \\
	y_1-y_2 = 0, \ y \geq \zeros
	\end{array} \right\} \\
    & = \left\{ \begin{bmatrix}
        0 & 0 & 0 \\
	0 & y_3-y_4 & 0 \\
	0 & 0 & -y_5
	\end{bmatrix} \middle| y_3-y_4 \geq 0, \ -y_5 \geq 0, \ y \geq \zeros
	\right\} \\
    & = \left\{ \begin{bmatrix}
        0 & 0 & 0 \\
	0 & t & 0 \\
	0 & 0 & 0
	\end{bmatrix} \middle| t \geq 0 \right\}.
    \end{align*}
Therefore, the matrix	
$$W = \begin{bmatrix}
		0 & 0 & 0\\
		0 & 1 & 0\\
		0 & 0 & 0
	\end{bmatrix}$$
    is a maximum-rank exposing vector, consistent with the description of $W^*$ in Lemma \ref{PFRD}. Let
	$$V = \begin{bmatrix}
		1 & 0\\
		0 &0\\
		0 &1\\
	\end{bmatrix} \in \mathbb{R}^{3 \times 2}$$
    be a corresponding facial range vector. 
    Then, the matrix variable $Y \in \mathbb{S}^3_+$ can be represented as $VRV^T$, where $R = \begin{bmatrix}
        Y_{00} & Y_{02} \\ Y_{20} & Y_{22}
    \end{bmatrix}$. 
    
    Indeed, the fact that $x_{1} = 0$ for every $x \in P$ is reflected in the above process: $R$ is obtained from $Y$ by removing the row and column corresponding to index $i = 1$.

    The resulting facially reduced problem is $\tilde{L} \cap \mathbb{S}^2_+$, where
    \[
    \tilde{L} = \left\{R = \begin{bmatrix}
        Y_{00} & Y_{02} \\
        Y_{20} & Y_{22}
    \end{bmatrix} \in \mathbb{S}^2 \, \middle| \, 
    \begin{array}{c}
         Y_{00} = Y_{02}, \, Y_{02} \geq 0,   \\
         Y_{00} = 1, \, Y_{22} = Y_{02}
    \end{array}
    \right\}
    = \left\{\begin{bmatrix}
        1 & 1 \\
        1 & 1
    \end{bmatrix}  \right\}.
    \]

    Let $\tilde{P} = \{x_2 \in \R \mid x_2 = 1, \, x_2 \geq 0\}$ be the polyhedron obtained from $P$ by eliminating $x_1$, which is fixed at 0 in $P$. It is straightforward to verify that $\tilde{L} \cap \mathbb{S}^2_+$ is the Shor relaxation for $\tilde{F} = \tilde{P} \cap \{0,1\}$. Applying partial FR with $\Ds$ to this reduced problem yields 
	$$\SEP(\tilde{L} \cap \Ds) = \left\{\begin{bmatrix}
		0 & 0 \\
		0 & y_{1}-y_{2}-y_{5}
	\end{bmatrix} \middle| \, \begin{array}{c} y_1 - y_2 - y_5 \geq 0, \\ y_1 - y_2 = 0, \, y \geq \zeros \end{array} \right\} = \{\mathbf{O}\}.$$
    This shows that there is no reduction in the second iteration, and the algorithm terminates with the final matrix order $r_P = 2$. 
\end{ex}

\begin{ex} \label{eg_ra}
	We now apply affine FR to the same set $F$ as in Example \ref{eg:pFRD}. The affine hull of $P$ is given by
	$$\aff P = \{ x\in \R^{2} \mid x_{1} + x_{2} = 1, x_{1} = 0\}.$$
	This yields a rank-$2$ exposing vector
    $$W = \begin{bmatrix}
		1 & 0\\
		-1 & 1\\
		-1 & 0
	\end{bmatrix}\begin{bmatrix}
		1 & 0\\
		-1 & 1\\
		-1 & 0
	\end{bmatrix}^{T} = \begin{bmatrix}
		1 & -1 & -1\\
		-1 & 2 & 1\\
		-1 & 1 & 1\\
	\end{bmatrix}.$$
	An associated facial range vector is 
	$$V  = \begin{bmatrix}
		1 \\
		0\\
		1\\
	\end{bmatrix} \in \mathbb{R}^{3 \times 1},$$
    which implies that the reduced matrix order is $r_A = 1$.
    
\end{ex}

\subsubsection{Diagonally dominant matrices} \label{sec_ddm}
Except for the zeroth row and column, the rows and columns of the matrix variable in the Shor relaxation correspond to the variables in the original binary formulation. We show below that, with the natural facial range vector described in Remark~\ref{rem:partial_fr_representation}, partial FR with $\K=\DDs$ removes the rows and columns corresponding to variables fixed at zero or one in $P$.

In analogy to Lemma \ref{PFRDlem1}, we first simplify the expression for $\SEP(L \cap \K)$ when $\K = \DDs$. Recall that $|u| = (|u_1|, \ldots, |u_n|)$ represents the componentwise absolute value of $u$.

\begin{lem}\label{apeq}
It holds that
\[
    \SEP(L \cap \DDs) = \left\{ \begin{bmatrix}
	-b^{T}y & \frac{1}{2}(y^{T}A - z^{T})  \\
	\frac{1}{2}(A^{T}y - z) & \Diag(z)\\
	\end{bmatrix} \, \middle| \,\begin{array}{c} 
        -b^T y \geq \frac{1}{2} \ones^T |A^Ty-z| \\
        z \geq \frac{1}{2}|A^Ty-z| \\ 
        y \in \R^m_+, \, z \in \R^n 
        \end{array}\right\}.
\]
\end{lem}
\begin{proof}
By Lemma \ref{shorAUX}, we have
    \begin{equation*}
		\SEP(L \cap \DDs) = \left\{ \begin{bmatrix}
			-b^{T}y & \frac{1}{2}(y^{T}A - z^{T})  \\
			\frac{1}{2}(A^{T}y - z) & \Diag(z)\\
		\end{bmatrix} \in \DD \middle| \, y \in \R^m_+, z \in \R^n \right\}.
    \end{equation*}
\linelabel{dddef}Recall that $\mathcal{D}\mathcal{D}$ is the set of diagonally dominant matrices and it is defined in \eqref{def_dd}. The above condition that the matrix is diagonally dominant is equivalent to the inequalities
\begin{equation} \label{eqn:SEP_DD}
    -b^T y \geq \frac{1}{2} \ones^T |A^Ty-z| \ \text{ and } \
        z \geq \frac{1}{2}|A^Ty-z|,
\end{equation}
which completes the proof.
\end{proof}

Given a polyhedron $P = \{x \in \R^n \mid Ax \leq b\}$, we partition $\{1, \ldots, n\}$ into the following mutually exclusive index sets: 
\begin{align*}
    I_+ &:= \{ i \in \{1,\ldots,n\} \mid x_{i} = 0 \text{ for all } x \in P\}, \\
    I_- &:= \{ i \in \{1,\ldots,n\} \mid x_{i} = 1 \text{ for all } x \in P\}, \\
    I_0 &:= \{ i \in \{1,\ldots,n\} \mid x_i \in (0,1) \text{ for some } x \in P\}.
\end{align*}
The following lemma shows that the row and column of an exposing vector corresponding to an index $i \in I_0$ must be zero.

\begin{lem}\label{PFRDD2}
If $W \in \SEP(L \cap \DDs)$, then $W_{ij} = W_{ji} = 0$ for all $i \in I_0$ and $j \in \{0, 1, \ldots, n\}$. Moreover, if $I_- = I_+ = \emptyset$, then $\SEP(L \cap \DDs) = \{\textbf{O}\}$.
\end{lem}
\begin{proof}
    By Lemma \ref{apeq}, if $W \in \SEP(L \cap \DDs)$, then there exist $y \in \R^m_+$ and $z \in \R^n$ satisfying \eqref{eqn:SEP_DD}. Let $w:= A^Ty$. Since $z \geq \frac{1}{2}|w-z|$, we have $4z_i^2 \geq (w_i-z_i)^2$ for all $i$, which is equivalent to 
    \[
    (w_i-3z_i)(w_i+z_i) \leq 0.
    \]
    If $w_i < 0$, then $w_i - 3z_i < 0$, and we must have $w_i+z_i \geq 0$ and 
    \[|w_i-z_i| = z_i - w_i \geq -2w_i.\]
    For any $x \in P$,
    \begin{align*}
        -b^Ty & \geq \frac{1}{2} \ones^T | w-z| \\
        & = \sum_{i: w_i < 0} \frac{1}{2}|w_i-z_i| + \sum_{i: w_i \geq 0} \frac{1}{2}|w_i-z_i| \\
        & \geq \sum_{i: w_i < 0} \frac{1}{2}|w_i-z_i| \\
        & \geq \sum_{i: w_i < 0} -w_i \\
        & \geq \sum_{i: w_i < 0} -w_i x_i + \sum_{i: w_i = 0} -w_i x_i + \sum_{i: w_i > 0} -w_i x_i \\
        & = -y^T Ax \\
        & \geq -b^Ty.
    \end{align*}
    Here, the fourth inequality holds due to Assumption \ref{ass1}. Since the first and last expression in this chain are equal, all inequalities must hold at equality. In particular, equality in 
    \[
    \sum_{i: w_i < 0} -w_i \geq \sum_{i: w_i < 0} -w_i x_i + \sum_{i: w_i = 0} -w_i x_i + \sum_{i: w_i > 0} -w_i x_i
    \]
    implies
    \[
    x_i = \begin{cases}
        1, & \text{ if } w_i < 0 \\
        0, & \text{ if } w_i > 0
    \end{cases} \quad \text{for all } x \in P.
    \]
    Hence, if $x_i \in (0,1)$ for some $x \in P$, then $w_i = 0$. Moreover, equality in 
    \[
    \sum_{i: w_i < 0} \frac{1}{2}|w_i-z_i| + \sum_{i: w_i \geq 0} \frac{1}{2}|w_i-z_i| \geq \sum_{i: w_i < 0} \frac{1}{2}|w_i-z_i|
    \]
    implies that $w_i = z_i$ whenever $w_i \geq 0$. Therefore, for $i \in I_0$, we have $z_i = w_i = 0$, which implies $W_{ij}=W_{ji} = 0$ for all $j \in \{0, 1, \ldots, n\}$. Finally, if $I_- = I_+ = \emptyset$, then $I_0 = \{1, \ldots, n\}$, so every entry of $W$ is zero except possibly $W_{00} = -b^Ty$. However, since $z_i = w_i = 0$ for all $i \in \{1, \ldots, n\}$, every term in the chain of inequalities equals zero, implying $-b^Ty = 0$. Therefore, $W = \mathbf{O}$ is the only element in $\SEP(L \cap \DDs)$.
\end{proof}

In analogy to Lemma \ref{PFRD}, the following result provides an explicit form of the exposing vectors generated by partial FR with $\K = \DDs$. We also show that the algorithm can terminate after at most one FR step. The key observation is that, after a single FR step, the new matrix variable $R$ is obtained from the original matrix variable $Y$ by removing the rows and columns corresponding to the variables that are fixed at zero or one in $P$. 

Without loss of generality, we assume that
\[
I_- = \{1, \ldots, p\} \ \text{ and } \ I_+ = \{p+1, \ldots, p+q\}.
\]
When $I_- = \emptyset$, we set $p = 0$; similarly, when $I_+ = \emptyset$, we set $q = 0$.

\begin{lem}\label{PFRDD}
Define the block-diagonal matrix $W^{*} \in \Snop$ as follows:
 \begin{equation} \label{maxW}
	W^{*} = \left[  
\begin{array}{c c c c c c c}  
p & -\textbf{1}^{T}& \vline & & \vline&  \\  
-\textbf{1} & I_{p} & \vline &  & \vline & \\
\hline   
 &  & \vline & I_{q}  &\vline & \\
 \hline
&  & \vline &   &\vline & \textbf{O}_{n-p-q} \\
\end{array}
\right],
\end{equation}  
where $I_{p}$ and $I_{q}$ are identity matrices of order $p$ and $q$, respectively, and $\mathbf{O}_{n-p-q}$ is the zero matrix of order $n-p-q$.
  
It holds that
	\begin{enumerate}
		\item $W^{*} \in \SEP(L \cap \DDs)$, and $\rank W^{*} = p + q$.
		\item $W^{*}$ has the maximum rank among all matrices in $\SEP(L \cap \DDs)$. 
		\item When applied to $L \cap \Snop$, partial FR with $\K = \DDs$ can terminate after at most one FR step.
	\end{enumerate}
\end{lem}
\begin{proof}
\begin{enumerate}
    \item  For each $i \in I_-$, since $x_i = 1$ for all $x \in P$, LP strong duality implies
\[
\max\{-b^Ty \mid A^T y = -e_i, \, y \geq \zeros\} = \min\{x_i \mid Ax \leq b\} = 1.
\]
Therefore, there exists $y^i \in \R^m_+$ such that $A^Ty^i = -e_i$ and $b^Ty^i = -1$. Similarly, for $i \in I_+$, since $x_i = 0$ for all $x \in P$, by \eqref{pfr_eq2}, there exists $y^i \in \R^m_+$ such that $A^Ty^i = e_i$ and $b^Ty^i = 0$. Define $y := \sum_{i \in I_- \cup I_+} y^i$ and $w := A^T y$. Then, 
\[
    w_i =  \begin{cases}
        -1 & \text{if } i \in I_-\\
        1 & \text{if } i \in I_+\\
        0 & \text{otherwise}.
    \end{cases}
\]
and $b^Ty = -p$. Let $z := |w| \in \Rn$ be the componentwise absolute value of $w$. It is straightforward to verify that $-b^Ty = p = \frac{1}{2} \ones^T|w-z|$ and $z \geq \frac{1}{2}|w-z|$. Therefore, by Lemma \ref{apeq},
	\begin{equation*}
	W^* = \begin{bmatrix}
			-b^Ty & \frac{1}{2}(w - z)^{T}  \\
			\frac{1}{2}(w - z) & \Diag(z)\\
		\end{bmatrix} \in \SEP(L \cap \DDs).
	\end{equation*}

To prove that $\rank W^{*} = p + q$, it suffices to show 
	 $$\rank \begin{bmatrix}
	     p & -\textbf{1}^{T}\\
	     -\textbf{1} & I_{p}\\
	 \end{bmatrix} = p.$$
 The rank of this matrix is at least $p$ as it contains the identity matrix $I_{p}$ as a submatrix, while the rank is at most $p$ since the all-ones vector in $\R^{p+1}$ is in the null space of this matrix.

\item Let $W \in \SEP(L \cap \DDs)$ be arbitrary. By Lemma \ref{PFRDD2}, the $i$-th row and column of $W$ must be zero for $i = p+q+1,\ldots,n$ (with indices starting from 0). Therefore, $W$ must be of the form
$$W = \begin{bmatrix}
    \bar{W} & \textbf{O}\\
    \textbf{O} & \textbf{O}
\end{bmatrix}\in \mathbb{S}_{+}^{n+1}$$
for some $\bar{W} \in \mathbb{S}_{+}^{p+q+1}$. It follows that $\rank W \leq p + q + 1$. Suppose, for the sake of contradiction, that $\rank W = p+q+1$, i.e., $\bar{W}$ is positive definite. Then a corresponding facial range vector $V \in \R^{(n+1) \times (n-p-q)}$ is given by
$$V = \begin{bmatrix}
    \textbf{O}\\
   I_{n-p-q}\\
\end{bmatrix}.$$
Now consider the facially reduced set 
\[
\tilde{L} := \{R \in \mathbb{S}^{n-p-q} \mid VRV^T \in L\},
\]
where $R$ is the reduced matrix variable. Observe that the first row of $V$ is zero, and hence for any $R \in \mathbb{S}^{n-p-q}$, we have $(VRV^{T})_{00} = 0$. However, the original arrow constraint in the Shor relaxation requires $Y_{00} = (VRV^{T})_{00} = 1$. This implies $\tilde{L} = \emptyset$, contradicting the assumption that $P \neq \emptyset$ in Assumption~\ref{ass1}. Therefore, no matrix in $\SEP(L \cap \DDs)$ can have rank greater than $p+q$, and thus $W^*$ has the maximum rank among all matrices in $\SEP(L \cap \DDs)$.

\item If $\SEP(L \cap \DDs) = \{\mathbf{O}\}$, then partial FR with $\K = \DDs$ terminates immediately. Otherwise, the algorithm generates a maximum-rank exposing vector in $\SEP(L \cap \DDs)$. Without loss of generality, assume that $W^*$ is selected to be the exposing vector. Define $r:=n - p - q$. A facial range vector $V \in \R^{(n+1) \times (r+1)}$ associated with $W^{*}$ is given by
$$V = \begin{bmatrix}
    \ones_{p+1} & \textbf{O}\\
    \textbf{0}_{q} & \textbf{O}\\
    \textbf{0}_{r} & I_{r}\\
\end{bmatrix}.$$

Let $R \in \mathbb{S}^{r+1}$ be the reduced matrix variable in the facially reduced formulation. We now derive the constraints of the facially reduced formulation. Recall that the original polyhedron $P$ is defined by the linear system $Ax \leq b$, where the $i$-th constraint is given by $a_{i}^{T}x \leq b_{i}$. Let $\tilde{A} \in \R^{m\times r}$ be the submatrix of $A$ consisting of its last $r$ columns, and denote by $\tilde{a}_{i}^{T}$ the $i$-th row of $\tilde{A}$. Define the modified right-hand side $\tilde{b} \in \R^{m}$ by 
$\tilde{b}_{i} = b_{i}-\sum_{j\in I_-}(a_i)_j$, for $i = 1, \ldots, m$. Then the data matrix corresponding to the $i$-th constraint in \eqref{shorith} becomes
\begin{equation*}
    V^{T}\begin{bmatrix}
		-b_{i} &  \frac{1}{2}a_{i}^{T} \\
		\frac{1}{2}a_{i} & \textbf{O}
\end{bmatrix}V = 
\begin{bmatrix}
		-\tilde{b}_{i} &  \frac{1}{2}\tilde{a}_{i}^{T} \\
		\frac{1}{2}\tilde{a}_{i} & \textbf{O}
\end{bmatrix} \in \mathbb{S}^{r+1}.
\end{equation*}
Under this reduction, the arrow constraint \eqref{def_arrow} transforms to $\tilde{\arrow}(R)=\tilde{e}_{0}$ with $\tilde{\arrow}: \mathbb{S}^{r+1} \rightarrow \R^{r+1}$ given by
\begin{equation*}
		\tilde{\arrow}(R) := \begin{bmatrix}
			R_{00} \\
			R_{11} - \frac{1}{2}(R_{01} + R_{10}) \\
			\vdots \\
			R_{rr} - \frac{1}{2}(R_{0r} + R_{r0}) \\
		\end{bmatrix} \in \mathbb{R}^{r+1},
	\end{equation*}
where $\tilde{e}_{0}$ is the first standard basis vector in $\R^{r+1}$. Let $\tilde{L}$ be the set defined by 
\[\left\langle \begin{bmatrix}
    -\tilde{b}_{i} &  \frac{1}{2}\tilde{a}_{i}^{T} \\
    \frac{1}{2}\tilde{a}_{i} & \textbf{O}
\end{bmatrix}, R \right\rangle \leq 0 \ \text{ for } i = 1,\ldots,m
\]
and $\tilde{\arrow}(R)=\tilde{e}_{0}$. Then the facially reduced reformulation of the Shor relaxation $L \cap \Snop$ is $\tilde{L} \cap \mathbb{S}_{+}^{r+1}$. On the other hand, $\tilde{L} \cap \mathbb{S}_{+}^{r+1}$ is exactly the Shor relaxation of the binary set 
$$\tilde{F} := \{0,1\}^{r} \cap \tilde{P}, \text{ where } \tilde{P}:= \{ y \in \R^{r} \vs \tilde{A}y \leq \tilde{b}\}.$$
Note that $$y \in \tilde{F} \text{  if and only if } \begin{bmatrix}
    \textbf{1}\\\textbf{0}\\y
\end{bmatrix} \in F.$$
By construction, none of the variables in $\tilde{P}$ is fixed at zero or one. Therefore, Lemma \ref{PFRDD2} implies that $\SEP(\tilde{L} \cap \DDs) = \{\textbf{O}\}$, and partial FR with $\K = \DDs$ terminates after this step. 
\end{enumerate}
\end{proof}



\begin{cor} \label{lem_rarp2}
It holds that $r_A \leq r_{P}^{+} \leq r_P$.
\end{cor}
\begin{proof}
    The proof follows the same reasoning as in Corollary \ref{lem_rarp}. Define two affine sets:
    \[
    P':= \left\{x \in \R^n \middle| \, \begin{array}{c} x_i = 0 \text{ for } i \in I_+ \\ x_i = 1 \text{ for } i \in  I_- \end{array} \right\} \text{ and } P'':= \{x \in \R^n \mid x_i = 0 \text{ for } i \in I_+\}.
    \]
    By definition of the index sets $I_+$ and $I_-$, it holds that $P \subseteq P' \subseteq P''$. Therefore, $\dim(P) \leq \dim(P') \leq \dim(P'')$. By construction, $r_A = \dim(P) + 1$, $r_P^+ = n+1 - |I_+| - |I_-| = \dim(P') + 1$ and $r_P = n+1 - |I_+| = \dim(P'') + 1$. Therefore, $r_A \leq r_P^+ \leq r_P$, as claimed.
\end{proof}

Note that the proof of Lemma \ref{PFRDD} also demonstrates that a maximum-rank exposing vector for $L \cap \DDs$ can be computed by solving a single LP. Therefore, its computational complexity is comparable to that of affine FR.

Corollary~\ref{lem_rarp2} shows that, for the Shor relaxation \eqref{shorSDP2} associated with any polyhedron $P=\{x\in\R^{n}\mid Ax\leq b\}$, affine FR produces a matrix variable whose order is no larger than those produced by partial FR with $\Ds$ or $\DDs$. Thus, the comparison is an exact, formulation-specific comparison for the Shor relaxation studied in this subsection.

The comparison in this subsection is limited to the cones $\D$ and $\DD$ and to the Shor relaxation \eqref{shorSDP2}. For affine-slice SDP relaxations of binary programs containing the normalization and arrow equations, \cite{huXu2026face} gives a broader structural analysis of diagonal-, diagonally-dominant-, and scaled-diagonally-dominant partial FR; in that setting, additional lifted affine equations can produce equality-group reductions and multiple partial-FR steps. For factor-width-$k$ partial FR, \cite{huSmithXu2026polyhedral} develops an LP-based face-identification method using subset-indexed moment variables and proves a comparison at every iteration. Unlike affine FR here, that construction depends on the chosen SDP relaxation and the factor width $k$.

The equality-generated SDP--RLT relaxations studied in \cite{hu2026sharpRLT} contain the additional lifted equations $AX=bx^{T}$ and are analyzed through their singularity degree. Inequality formulations represented using nonnegative slacks instead lead to a product cone involving both a positive semidefinite cone and a nonnegative orthant. These works therefore address different formulations and questions from the one-step $\D$- and $\DD$-partial-FR comparison established here.


\subsection{Sieve-SDP}\label{sec_sieve}
In \cite{zhu2019sieve}, Zhu et al. propose a specialized FRA called \emph{Sieve-SDP}. Their method aims to efficiently identify certain structural patterns in the constraint matrices that allow for low-cost extraction of exposing vectors. The facial reduction component of Sieve-SDP can be summarized as follows: for each constraint of the form $\langle A_{i}, Y \rangle = b_{i}$ or $\langle A_{i}, Y \rangle \leq b_{i}$, the algorithm checks whether the matrix $A_{i}$ can be permuted into the block-diagonal form:
\begin{equation}\label{sieve_sdp_eq1}
	A_{i} = \begin{bmatrix}
		D_{i} & \textbf{O} \\
		\textbf{O} & \textbf{O}
	\end{bmatrix} \text{ with } D_{i} \in \mathbb{S}_{++}^{k}.
\end{equation}
If $b_{i} = 0$, then $A_{i}$ is an exposing vector, and facial reduction can be applied by eliminating the rows and columns of the matrix variable $Y$ corresponding to the block $D_{i}$. The algorithm terminates when no such constraint can be found. Sieve-SDP effectively searches for specific elements within the set of exposing vectors and relies only on an incomplete Cholesky factorization to verify the positive definiteness of $D_{i}$. Each such test is therefore computationally inexpensive.

Next, we analyze the behavior of Sieve-SDP when applied to the Shor relaxation.
\begin{lem}\label{lem_sieve}
    Let $F = P \cap \{0,1\}^{n}$, where $P = \{x \in \R^n \mid Ax \leq b\}$, and let $L \cap \Snop$ denote the Shor relaxation for $F$, where $L$ is defined in \eqref{shorSDP2}. Under Assumption \ref{ass1}, Sieve-SDP does not perform any reduction when applied to $L \cap \Snop$. In particular, $r_{S} = n+1$.
\end{lem}
\begin{proof}
	We examine all data matrices in the Shor relaxation \eqref{shorSDP2} to determine whether Sieve-SDP can identify any exposing vectors. Note that a matrix can be permuted into the form \eqref{sieve_sdp_eq1} only if it is positive semidefinite.
	The data matrices corresponding to the linear constraints $\mathcal{A}(Y) \leq \textbf{0}$ are
	\begin{equation}\label{sieve1}
		A_i:= \begin{bmatrix}
			-b_{i} &  \frac{1}{2}a_{i}^{T} \\
			\frac{1}{2}a_{i} & \textbf{O}
		\end{bmatrix} \text{ for } i =1,\ldots,m.
	\end{equation}
    For any $A_i$ to be positive semidefinite, it must be that $a_{i} = 0$. However, since $P \neq \emptyset$ by Assumption \ref{ass1}, this implies $b_{i} \geq 0$. Consequently, $A_i \notin \Snop \setminus \{\mathbf{O}\}$, so no matrix in \eqref{sieve1} yields an exposing vector.
	
    Next, consider the data matrices associated with the arrow constraint $\text{arrow}(Y) = e_{0}$:
	\begin{equation}\label{sieve2}
		\begin{bmatrix}
			1 & \textbf{0}^{T}\\
			\textbf{0} & \textbf{O}
		\end{bmatrix} \text{ and } 	\begin{bmatrix}
			0 & -\frac{1}{2}e_{i}^{T}\\
			-\frac{1}{2}e_{i} & e_{i}e_{i}^{T}
		\end{bmatrix} \text{ for } i =1,\ldots,n.
	\end{equation}
    The first matrix is positive semidefinite, but its corresponding right-hand side is 1, and Sieve-SDP requires a right-hand side of zero to eliminate rows and columns. The remaining matrices in \eqref{sieve2} are not positive semidefinite. 
    
    Therefore, none of the constraints satisfy the criteria for exposing vectors under Sieve-SDP, and no reduction occurs. It follows that $r_S = n+1$.
\end{proof}

Combining Corollary~\ref{lem_rarp2} with Lemma~\ref{lem_sieve} proves Theorem~\ref{thm_main}.

One possible extension of Sieve-SDP is to use valid equality constraints derived from the original problem. Specifically, for an equality constraint $a_{i}^{T}x = b_{i}$ ($a_i \neq \zeros$) in the original problem, one may add the valid equality constraint $\langle A_i,Y \rangle = 0$, where 
$$A_i := \begin{bmatrix}
	-b_{i}\\
	a_{i}
\end{bmatrix}\begin{bmatrix}
	-b_{i}\\
	a_{i}
\end{bmatrix}^{T}.$$
Since $A_i$ is nonzero and positive semidefinite, it is an exposing vector for the strengthened relaxation. A modified procedure could aggregate the matrices obtained from several such equalities and compute a facial range vector from the null space of their sum. Unlike standard Sieve-SDP, however, a singular positive semidefinite matrix does not in general justify deleting its entire supported principal block; the reduction must use its null space. Investigating the effectiveness and computational performance of this extension is an interesting direction for future research.

\section{Numerical Experiments}\label{sec_bip}
In this section, we evaluate the performance of affine FR and compare it with partial FR. The computational study is designed to separate two related but distinct questions:
\begin{enumerate}[label=(\roman*), itemsep=0pt, topsep=2pt]
	\item whether affine FR can be computed efficiently and produces substantial reductions in the order of the matrix variable across a broad collection of binary and mixed-binary formulations;
	\item whether, on instances and problem classes for which the resulting SDP relaxations are within reach of current solvers, preprocessing affects SDP solution time and numerical reliability.
\end{enumerate}
The MIPLIB experiments address the first question, while the weighted SAT experiments address the second on instances for which the SDP relaxation can be substantially stronger than the LP relaxation.

All experiments were conducted on a Mac Studio equipped with an Apple M2 Ultra processor, 128 GB of RAM, and macOS. We use Gurobi to solve the LP problems that arise in preprocessing. For solving the original SDP relaxations and the relaxations generated by affine FR, we use Mosek (version 11.0) \cite{mosek} in conjunction with YALMIP \cite{lofberg2004yalmip}. For affine FR, we compute the affine hull by solving the LP in \eqref{splitFR}. When constructing a facial range vector from an eigendecomposition, eigenvalues below \(10^{-5}\) are treated as zero. This conservative threshold may reduce the amount of reduction in borderline cases, but it is less likely to remove feasible directions. We impose a 5-minute time limit for all LPs used in preprocessing. If the solver fails to return a solution within this time limit, we classify the instance as yielding no reduction.

\subsection{MIPLIB mixed-binary experiments}\label{sec_miplib}
This subsection evaluates the preprocessing cost and matrix-order reductions obtained by affine FR on a broad collection of MIPLIB mixed-binary formulations. We study the preprocessing step rather than the strength or end-to-end solution of the Shor relaxation, which is often beyond the reach of generic SDP solvers at this scale.

We consider mixed-binary linear programming instances from the MIPLIB 2017 \cite{MIPLIB2017}, specifically those under the ``collection'' category with at most $10,000$ variables. The instances are read into MATLAB R2024b using Gurobi (version 11.0.0) \cite{gurobi}. 
We include only instances whose discrete variables are binary; integer variables bounded between zero and one are treated as binary variables.
We also exclude seven instances due to numerical issues\footnote{For these instances, the LP solver reported numerical difficulties when solving \eqref{splitFR}. The excluded instances are \texttt{control20-5-10-5}, \texttt{control30-5-10-4}, \texttt{app1-1}, \texttt{app2-1}, \texttt{app2-2}, \texttt{lrn}, and \texttt{rentacar}.} encountered when computing the affine hull of $P$.

The final dataset consists of $332$ mixed-binary linear programming instances, drawn from multiple application areas and spanning a range of problem sizes for testing SDP preprocessing. The MIPLIB 2017 instances are available at \url{https://miplib.zib.de/}. Our source code, including the list of tested instances, can be found at \url{https://github.com/haohu-code/affineFR_code}.

The mixed-binary linear programming instances from \cite{MIPLIB2017} take the following form:
\begin{equation}\label{mblp}
	\begin{array}{rllr}
		\min & c^{T}x &\\
		\text{subject to} & Ax = b &\\
		&Bx \leq d& \\
		& l \leq x \leq u\\
		& x_{i} \in \{0,1\}  &\text{ for } i \in \mathcal{B},
	\end{array}
\end{equation}
where $x \in \Rn$ is the decision variable. The linear system $Ax = b$ and $Bx \leq d$ represent the equality and inequality constraints, respectively.  The vectors $l$ and $u$ specify the lower and upper bounds for $x$, respectively. If $l_{i} = -\infty$, then $x_{i}$ has no lower bound. Similarly, $u_{i} = \infty$ means $x_{i}$ has no upper bound. The subset $\mathcal{B} \subseteq \{1,\ldots,n\}$ indexes the binary variables.

To apply partial FR, we use the mixed-binary version of the Shor relaxation. Let $L$ be defined by

\begin{equation}\label{shorSDP}
L:=\left\{Y=\begin{bmatrix}
			1 & x^{T} \\
			x & X
		\end{bmatrix}\in\Sno \,\middle|\,
        \begin{array}{l}
        Ax=b,\; Bx\leq d,\; l\leq x\leq u,\\
        x_i=X_{ii}\ \text{for every } i\in\mathcal B
        \end{array}
        \right\}.
\end{equation}

The corresponding SDP relaxation is $L\cap\Snop$, and the order of the matrix variable $Y$ is $n+1$. This formulation specifies how partial FR is applied; affine FR itself is independent of this particular SDP relaxation and can be applied to any SDP relaxation to produce a no-weaker relaxation with a reduced matrix variable. Since affine FR and the two partial-FR variants considered below are all LP-based preprocessing procedures with comparable theoretical complexity, the comparison focuses on the reductions in matrix order obtained and the preprocessing time required.

To apply partial FR in the MIPLIB experiments, we implement the algorithm described in \cite{permenter2018partial} to identify maximum-rank exposing vectors from $\SEP(L \cap \Ds)$ and $\SEP(L \cap \DDs)$.\footnote{The original code in \cite{permenter2018partial} targets generic SDP relaxations. For the cone of diagonally dominant matrices $\DD$, a direct implementation introduces $\mathcal{O}(n^{2})$ additional constraints. By exploiting the arrowhead structure in the Shor relaxation, our customized implementation reduces this to $\mathcal{O}(n)$ constraints.}
For each FRA, we construct a facial range vector $V \in \mathbb{R}^{(n+1) \times r}$ with $r \leq n + 1$, and replace $Y$ with a smaller matrix variable $R \in \mathbb{S}_{+}^{r}$ via the substitution $Y= VRV^{T}$ (see \eqref{K1fr}). As shown in the results below, the order of the matrix variable can be reduced through facial reduction.

We next summarize the reductions in matrix order and preprocessing times.\phantomsection\label{sec_ms} The aggregate results are reported in Tables \ref{froverview} and \ref{frtime}. Complete per-instance results, including solver statuses and separate times for the auxiliary LP and facial-range-vector construction, are provided with the source code in machine-readable format.

Table \ref{froverview} summarizes the performance of the three FRAs on all 332 instances.
\begin{itemize}
    \item ``Reduced instances'' reports the number and percentage of instances for which the reduced order satisfies $r<n+1$.
    \item \linelabel{AFRed_l} \label{AFRed_p}``Mean ratio (all)'' reports the mean of $r/(n+1)$ across all 332 instances, while ``Median ratio (reduced)'' reports the median among instances with $r<n+1$. A smaller ratio indicates a greater reduction.
    \item ``Time limits'' and ``Other failures'' report unsuccessful preprocessing runs by status.
\end{itemize}
As shown in Table \ref{froverview}, affine FR reduces substantially more instances than either partial-FR variant and gives a smaller mean order ratio, indicating larger reductions in matrix order on this test set. All three methods complete every preprocessing run without a time limit or other failure.

\begin{table}[htbp]
	\centering
	\small
	\begin{tabular}{@{}lrrrrr@{}}
	\toprule
	Method & \shortstack{Reduced\\instances} & \shortstack{Mean ratio\\(all)} & \shortstack{Median ratio\\(reduced)} & \shortstack{Time\\limits} & \shortstack{Other\\failures} \\
	\midrule
	Affine FR & 265 (79.8\%) & 0.8577 & 0.8688 & 0 & 0 \\
	Partial FR ($\Ds$) & 42 (12.7\%) & 0.9862 & 0.9069 & 0 & 0 \\
	Partial FR ($\DDs$) & 52 (15.7\%) & 0.9815 & 0.9200 & 0 & 0 \\
	\bottomrule
	\end{tabular}
	\caption{Matrix-order reductions across the 332 MIPLIB mixed-binary instances. Ratios are reduced order divided by original order.}\label{froverview}
\end{table}

Table \ref{frtime} reports the median and mean total preprocessing time for each FRA, grouped by the number of variables $n$. Each entry is displayed as median (mean) in seconds and includes both the computation of an exposing vector and the construction of the corresponding facial range vector. The number of instances in each size group is shown explicitly. Affine FR has the smallest median and mean time in every size group in this run.

In particular, the MIPLIB results concern the reduction step itself; they are not intended to show that the unreduced or reduced Shor relaxation is a competitive general-purpose method for solving MIPLIB instances.

\begin{table}[htbp]
	\centering
	\small
	\begin{tabular}{lrrrr}
	\toprule
	Number of variables & Instances & Affine FR & Partial FR ($\Ds$) & Partial FR ($\DDs$) \\
	\midrule
	$n\leq 2000$ & 140 & 0.05 (0.09) & 0.08 (0.19) & 0.08 (0.19) \\
	$2000<n\leq 4000$ & 84 & 0.69 (0.93) & 1.53 (1.70) & 1.55 (2.23) \\
	$4000<n\leq 6000$ & 54 & 2.11 (2.20) & 4.07 (4.42) & 4.07 (4.45) \\
	$6000<n\leq 8000$ & 27 & 5.62 (5.39) & 7.91 (8.67) & 8.08 (8.87) \\
	$8000<n\leq 10000$ & 27 & 11.82 (10.69) & 16.19 (15.97) & 16.24 (19.42) \\
	\bottomrule
	\end{tabular}
	\caption{Preprocessing time in seconds by problem size. Each entry is median (mean) total preprocessing time.}\label{frtime}
\end{table}

\subsection{SAT and weighted SAT experiments: SDP bounds and solver reliability}\label{sec_sat}
This subsection reports the SAT part of the numerical study. Subsection~\ref{appB3} reports the screening time and matrix-order reductions across the SATLIB instances, while Subsection~\ref{appB4} examines SDP solution time and solver behavior on selected weighted-SAT instances with substantial reductions.

\subsubsection{The weighted SAT problem}\label{appB1}
The SAT problem is a fundamental problem in combinatorial optimization and computer science. While the standard SAT problem simply asks whether there exists an assignment of truth values to a set of Boolean variables that satisfies a given propositional formula (typically expressed in Conjunctive Normal Form), the \textit{weighted} SAT problem additionally seeks to minimize a specified cost associated with that satisfying assignment. 

Let $\mathcal{V} = \{1, \dots, n\}$ be the set of Boolean variables, and let $\mathcal{C} = \{C_1, \dots, C_m\}$ be the set of clauses. For each variable $i \in \mathcal{V}$, let $w_i \in \mathbb{R}$ denote its associated weight. For each clause $C_j \in \mathcal{C}$, let $C_j^+$ denote the set of indices of variables appearing positively, and let $C_j^-$ denote the set of indices of variables appearing negatively. We introduce binary decision variables $x_i \in \{0, 1\}$ for each $i \in \mathcal{V}$, where $x_i = 1$ if the variable is true and $x_i = 0$ if it is false. 

To satisfy the formula, every clause must contain at least one true literal. Combining this logical satisfaction requirement with the objective function, the minimization version of the weighted SAT problem is formulated as the following binary integer linear program:
\begin{equation} \label{eq:weighted_sat}
\begin{array}{rll}
\text{Minimize} & \displaystyle \sum_{i=1}^{n} w_i x_i \\
\text{subject to} & \displaystyle \sum_{i \in C_j^+} x_i + \sum_{i \in C_j^-} (1 - x_i) \geq 1, & \text{for } j = 1, \dots, m, \\
& x_i \in \{0, 1\}, & \text{for } i = 1, \dots, n.
\end{array}
\end{equation}

\subsubsection{The LP and SDP relaxations}\label{appB2}

The continuous LP relaxation of the weighted SAT problem is obtained by replacing the binary constraints $x_i \in \{0, 1\}$ with the continuous interval $x_i \in [0, 1]$. While modern commercial solvers can process this linear relaxation highly efficiently, the resulting bounds can be substantially weaker than SDP bounds on some weighted SAT instances.

This difference in bound strength motivates the use of SDP relaxations, which can, in principle, deliver much tighter bounds than the corresponding linear relaxation. Substantial evidence for the effectiveness of SDP-based methods for SAT and related Boolean problems appears in \cite{anjos2005improved,de2000relaxations,sinjorgo2024solving,gomes2006power}. In addition, one can strengthen an SDP relaxation by adding valid constraints, as in the RLT formulation below. The key question is whether the resulting SDP can actually be solved in practice. That is precisely where the main difficulties arise, including large order of the matrix variable and numerical difficulties associated with the failure of Slater's condition. The experiments in Subsection~\ref{appB4} compare the solution time and solver status of the original and reduced formulations on the selected weighted SAT instances.

To construct a stronger SDP relaxation, we strengthen the Shor relaxation \eqref{shorSDP2} using the Reformulation-Linearization Technique (RLT). Let $Q \in \mathbb{S}^{n+1}$ be the cost matrix defined such that its first row and column correspond to $\frac{1}{2}w$ (with $Q_{00} = 0$ and all other entries being zero), ensuring that the objective evaluates to $\langle Q, Y \rangle = \sum_{i=1}^{n} w_i x_i$.

For the $j$-th linear inequality $a_{j}^{T}x \geq b_{j}$ in \eqref{eq:weighted_sat}, we define the extended coefficient vector $\hat{a}_j = \begin{bsmallmatrix} -b_{j} \\ a_{j} \end{bsmallmatrix} \in \mathbb{R}^{n+1}$. To capture the Boolean nature of the variables more tightly, we impose additional valid inequalities on the Shor relaxation:
\begin{enumerate}
    \item \textbf{Element-wise nonnegativity:} Since $x_i \in \{0, 1\}$, any product of these variables is non-negative, which is enforced by $Y \geq \mathbf{O}$ element-wise.
	\item \textbf{RLT-type inequalities:} Multiplying both sides of the original linear constraints by the non-negative vector $x$ yields valid inequalities. In the lifted matrix space, these correspond to $Y\hat{a}_{j} \geq \zeros$ for all $1 \leq j \leq m$.
\end{enumerate}

Combining these ingredients gives the following RLT-strengthened SDP relaxation, obtained by starting with the Shor relaxation and adding the above inequalities:
\begin{equation} \label{eq:sat_sdp_rlt}
\begin{array}{rll}
\min & \langle Q, Y \rangle \\
\text{subject to} & \arrow(Y) = e_{0} \\
& Y\hat{a}_{j} \geq \zeros & \text{for } 1 \leq j \leq m \\
& Y \in \Snop \\
& Y \geq \mathbf{O}.
\end{array}
\end{equation}

The element-wise non-negativity constraints $Y \geq \mathbf{O}$ introduce a significant computational burden for interior-point SDP solvers. We retain them because they strengthen the relaxation, despite this additional computational cost.

\subsubsection{Affine FR applied to the SAT problem}\label{appB3}
When affine FR is applied to this formulation, the matrix variable $Y\in \Snop$ is replaced by the lower-dimensional substitution $Y=V R V^T$ where $R \in \Srp$, effectively reducing the order of the matrix variable while preserving a valid, no-weaker RLT-type SDP relaxation.

For the large-scale experiment in this subsection, we do not form an SDP or construct the facial range matrix $V$. After solving the auxiliary LP, we determine the reduced matrix order directly from the numerical rank of the matrix of affine equations, computed by sparse QR factorization. The reported computation time includes the auxiliary LP and this sparse-rank calculation. The facial range matrix is constructed only for the weighted SAT instances in Subsection~\ref{appB4}, where reduced SDP relaxations are solved.

We downloaded all 97 benchmark archives available from SATLIB~\cite{satlib}, comprising 50,306 SAT instances. We parsed each instance in DIMACS format, retained instances having at most 100,000 variables and 150,000 clauses, and applied the split auxiliary LP used by affine FR to every remaining parseable instance, without solving an SDP. We retained a family when at least one of its instances produced a nonzero exposing matrix. All size-eligible instances from each retained family, including those with no reduction, were then included in the family-level aggregates. This procedure selected the ten families \texttt{Bejing}, \texttt{QG}, \texttt{bf}, \texttt{blocksworld}, \texttt{bmc}, \texttt{hanoi}, \texttt{jnh}, \texttt{logistics}, \texttt{parity}, and \texttt{ssa}.

\begin{table}[htbp]
\centering
\resizebox{\textwidth}{!}{%
\begin{tabular}{lcccccc} \hline
\multirow{2}{*}{Data Set Name} & \multirow{2}{*}{\# Instances} & \multirow{2}{*}{\# Reduced} & \multicolumn{2}{c}{Average Matrix Order} & \multirow{2}{*}{Avg Ratio} & \multirow{2}{*}{\shortstack{Avg LP+QR\\Time (s)}} \\ \cline{4-5} & & & Original ($n+1$) & Reduced ($r$) & & \\ \hline
QG & 22 & 22 & 1044 & 592 &  0.55 &  0.45 \\ 
bf & 4 & 4 & 1699 & 644 &  0.39 &  0.23 \\ 
blocksworld & 7 & 7 & 1645 & 1620 &  0.94 &  0.19 \\ 
bmc & 7 & 6 & 12416 & 8627 &  0.56 & 49.06 \\ 
hanoi & 2 & 2 & 1326 & 1013 &  0.76 &  0.08 \\ 
Bejing & 16 & 6 & 9120 & 7281 & 0.91 & 3.25 \\ 
jnh & 50 & 27 & 101 & 99 &  0.98 &  0.02 \\ 
logistics & 4 & 4 & 1882 & 1783 &  0.94 &  0.10 \\ 
parity & 30 & 15 & 1045 & 576 &  0.66 &  0.16 \\ 
ssa & 8 & 8 & 2352 & 277 &  0.23 & 28.34 \\ 
\hline
\end{tabular}
}
\caption{Average matrix-order reductions and LP-plus-rank computation times for SAT datasets}\label{tab:sat_results_2}
\end{table}

Table \ref{tab:sat_results_2} reports the aggregate preprocessing results achieved by affine FR for these ten families. The columns report the following aggregated metrics:
\begin{itemize}
    \item ``Data Set Name'': The category or family of SAT instances.
    \item ``\# Instances'': The number of instances evaluated in the family.
    \item ``\# Reduced'': The number of instances for which affine FR identified an exposing vector and reduced the order of the matrix variable.
    \item ``Original ($n+1$)'': The average original order of the SDP matrix variable (rounded to the nearest integer) across instances in the family; this equals the number of Boolean variables plus one.
    \item ``Reduced ($r$)'': The average reduced order of the matrix variable (rounded to the nearest integer) after applying affine FR.
    \item ``Avg Ratio'': The mean reduction ratio, defined as the reduced matrix order ($r$) divided by the original matrix order ($n+1$). A lower ratio indicates a larger reduction in the order of the matrix variable.
    \item ``Avg LP+QR Time (s)'': The average computation time, in seconds, required to solve the auxiliary LP and compute the reduced matrix order by sparse QR factorization. It excludes construction of an explicit facial range matrix.
\end{itemize}

The auxiliary LP reached optimality for 50,296 of the 50,306 downloaded instances; eight instances exceeded the size limits, one was nonconforming, and one reached the time limit. Of those 50,296 instances, affine FR identified nonzero exposing matrices in 101, while the remaining 50,195 yielded no detected reduction under the studied SAT formulation. All 101 reductions occurred among the 150 instances in the ten families reported in Table~\ref{tab:sat_results_2}. The existence and rank of an exposing vector depend on the affine structure of the formulation.

The average LP-plus-rank time is below one second for seven of the ten families. The \texttt{Bejing} family averages 3.25 seconds. The larger averages for \texttt{bmc} and \texttt{ssa} are driven by individual difficult instances: \texttt{bmc-ibm-4}, with 28,161 variables and 139,716 clauses, reaches the five-minute time limit, while \texttt{ssa6288-047}, with 10,410 variables and 34,238 clauses, requires 226.26 seconds. These results show that the sparse screening computation remains modest for most instances, while the auxiliary LP can be costly on particular large formulations.

\subsubsection{Weighted SAT bounds and solver behavior}\label{appB4}
Table~\ref{tab:sat_results_3} reports computational results for the minimization version of the weighted SAT problem. From the 150 instances summarized in Table~\ref{tab:sat_results_2}, the comparison includes every instance with at most 1,750 variables and 2,000 clauses for which the reduced matrix order satisfies $r<0.95(n+1)$ and $r\leq350$. These criteria depend only on the instance size and the reduction obtained by affine FR; they do not use the objective coefficients or the resulting LP and SDP bounds. Fourteen instances satisfy these criteria. In these experiments, we generated the weights as $n$ independent nonnegative integer values, each drawn as the absolute value of a standard normal random variable, scaled by a factor of 10, and then rounded to the nearest integer. The resulting fixed objective vectors are provided with the source code and are used in all reported runs.

The table compares the bounds obtained from the LP relaxation of \eqref{eq:weighted_sat}, the SDP relaxation \eqref{eq:sat_sdp_rlt}, and the relaxation generated by affine FR from \eqref{eq:sat_sdp_rlt}.

For each test instance, the table reports the true optimal value (or ``Inf'' if the instance is infeasible). It then compares the optimal value with the continuous LP bound. For both the original SDP relaxation and the relaxation generated by affine FR, we report the computed lower bound and the order of the matrix variable. A bold number indicates that the relaxation matched the known optimal value to displayed precision. The label ``OOM'' indicates that a preliminary attempt to construct and solve the unreduced SDP terminated the MATLAB process without returning a solver result, consistent with memory exhaustion. After observing these failures, we imposed an order limit of 350 on unreduced SDP models in subsequent batch runs and executed the SDP solves in separate MATLAB processes.

The numerical results highlight four points about applying affine FR to SDP relaxations of binary optimization problems:

\begin{itemize}
    \item \textbf{Strength of the bound:} On these weighted SAT instances, the SDP relaxation often provides a much tighter lower bound than the LP relaxation when it can be solved reliably. For smaller feasible instances such as \texttt{anomaly} and \texttt{medium}, the SDP bounds recover the exact optimal value.
      
    \item \textbf{Overcoming memory bottlenecks:} For the \texttt{par8} family and \texttt{ssa0432-003}, which have original matrix orders of $351$ and $436$, respectively, preliminary attempts to construct and solve the unreduced SDP terminated MATLAB without returning a solver result. After applying affine FR, these orders are substantially reduced (e.g., from $351$ to $65$), and the reduced SDP models are solved in a matter of seconds.

	\item \textbf{Solution time:} For the displayed \texttt{jnh} instances, solving the relaxation generated by affine FR is faster than solving the original SDP in the reported run, although the improvement varies by instance.
    
   \item \textbf{Numerical reliability:} 
   Because the unreduced SDP may fail Slater's condition, differences between the reported solver outcomes before and after applying affine FR should be interpreted with care. They can range from small changes in reported bounds to different solver statuses, such as a finite reported value versus infeasibility. Such differences may reflect numerical regularity effects and, when the affine-FR face restriction is not already implied by the original relaxation, the strengthening introduced by the relaxation generated by affine FR; see \cite{permenter2017reduction,sremac2021error,tunccel2016polyhedral}. This behavior is illustrated by the infeasible instance \texttt{jnh309}. When solving its original unregularized SDP, Mosek reports numerical issues and returns a finite lower bound; for the relaxation generated by affine FR, Mosek instead reports infeasibility without a numerical warning.

\end{itemize}

Overall, these experiments show that the SDP relaxation can be substantially stronger than the LP relaxation on the weighted SAT instances tested here. On several tested instances, preliminary attempts to solve the unreduced SDP terminated MATLAB or the solver produced numerical warnings. Affine FR reduces the order of the matrix variable and, in these experiments, improves the observed solution time and solver behavior.

\begin{table}[htbp]
\scriptsize
\centering
\resizebox{\textwidth}{!}{%
\begin{tabular}{lc|c|ccc|ccc} \hline
\multirow{2}{*}{Instance} & \multirow{2}{*}{Optimal} & \multirow{2}{*}{LP bound} & \multicolumn{3}{c|}{Original SDP} & \multicolumn{3}{c}{Affine FR} \\ 
&  &  & SDP bound & Matrix Order & Time & SDP bound & Matrix Order & Time \\ \hline
        anomaly &   174.00 &   122.50 & \textbf{174.00}&    49 &    0.20 &  \textbf{174.00}&    39 &    0.10 \\ 
         medium &   225.00 &   188.75 & \textbf{225.00}&   117 &    7.02 &  \textbf{225.00}&   107 &    4.43 \\ 
          jnh14 &      Inf &   261.92 & 433.41   &   101 &    9.88 &  435.73   &    94 &    6.79 \\ 
         jnh210 &   302.00 &   206.50 & 221.85   &   101 &    9.32 &  221.85   &    95 &    5.64 \\ 
         jnh215 &      Inf &   182.36 & 280.62   &   101 &   23.93 &  280.64   &    95 &   14.14 \\ 
         jnh302 &      Inf &   257.75 &    \textbf{Inf}&   101 &    6.71 &     \textbf{Inf}&    93 &    4.23 \\ 
         jnh309 &      Inf &   262.00 & 414.14   &   101 &    7.74 &     \textbf{Inf}&    95 &    6.53 \\ 
           jnh7 &   334.00 &   220.74 & 258.67   &   101 &   39.84 &  258.18   &    95 &   12.06 \\ 
         par8-1 &   666.00 &   388.00 &    OOM   &   351 &       - &  539.29   &    65 &    1.83 \\ 
         par8-2 &   611.00 &   395.50 &    OOM   &   351 &       - &  581.51   &    69 &    1.43 \\ 
         par8-3 &   513.00 &   278.33 &    OOM   &   351 &       - &  426.73   &    76 &    2.29 \\ 
         par8-4 &   677.00 &   339.70 &    OOM   &   351 &       - &  499.50   &    68 &    1.41 \\ 
         par8-5 &   609.00 &   331.42 &    OOM   &   351 &       - &  530.71   &    76 &    2.88 \\ 
    ssa0432-003 &      Inf &  1179.29 &    OOM   &   436 &       - & 1284.79   &   117 &   16.04 \\ 
\hline
\end{tabular}
}
\caption{LP and SDP results before and after affine FR}\label{tab:sat_results_3}
\end{table}

\section{Conclusion}\label{sec_con}
We introduced affine FR, an automatic, LP-based preprocessing method for SDP relaxations of binary and mixed-binary optimization problems. The method uses $\aff P$ to construct a face of the positive semidefinite cone containing the lifted feasible set and applies the substitution $Y=VRV^T$ to reduce the order of the matrix variable. Because the constructed face need not contain the entire feasible region of the original SDP relaxation, the relaxation generated by affine FR may strengthen the original relaxation. When $\aff P=\aff F$, the relaxation generated by affine FR satisfies Slater's condition.

For the Shor relaxation of a binary program, the theoretical comparison establishes that affine FR produces a matrix variable whose order is no larger than those produced by the two partial-FR variants and Sieve-SDP considered here. On MIPLIB 2017 mixed-binary instances, affine FR reduces matrix order more often and by a larger amount than the partial-FR variants considered, with comparable preprocessing times. On selected SAT relaxations with substantial reductions, affine FR shortens SDP solution times and yields more reliable solver outcomes when the unreduced formulations encounter numerical difficulties. These results concern the preprocessing of SDP relaxations; the usefulness of the resulting formulation also depends on the strength of the chosen relaxation and the numerical behavior of the solver.
\section*{Acknowledgements}
The authors would like to thank Professor G\'abor Pataki for his valuable insights and comments on Sieve-SDP, which have contributed to the improved presentation of this work. ChatGPT was utilized to assist with language polishing of this manuscript. Its use was limited to improving linguistic clarity and did not contribute to the development of the manuscript's scientific content.

\section*{Declarations}

\subsection*{Funding}
This work was supported by the Air Force Office of Scientific Research under award number FA9550-23-1-0508.

\subsection*{Competing interests}
The authors declare that they have no competing interests.

\subsection*{Data availability}
The MIPLIB 2017 instances used in the numerical experiments are publicly available at \url{https://miplib.zib.de/}. Complete per-instance MIPLIB results are provided with the source code as a CSV file. The SAT benchmark instances used in Subsection~\ref{sec_sat} are publicly available from SATLIB~\cite{satlib}. The fixed weight vectors for the weighted SAT experiments, together with the source code and tested-instance lists needed to reproduce the experiments, are available at \url{https://github.com/haohu-code/affineFR_code}.

\subsection*{Code availability}
The source code for the computational experiments, including the list of tested instances, is available at \url{https://github.com/haohu-code/affineFR_code}.

\bibliographystyle{unsrtnat}
\bibliography{mybib}

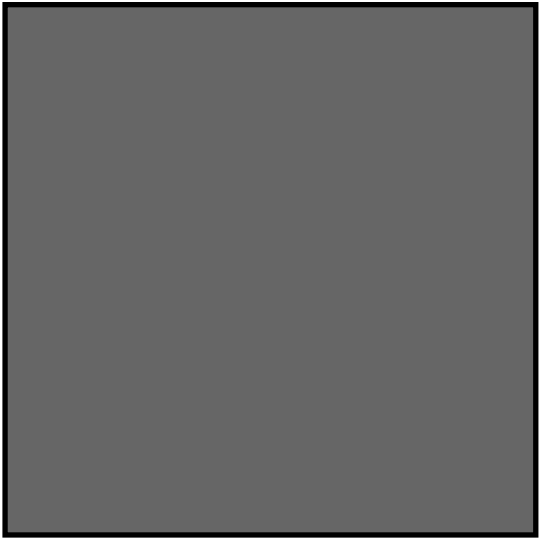
\end{document}